\documentclass[10pt]{article}

\usepackage{amsthm}
\usepackage{amsmath}
\usepackage{amsfonts}
\usepackage{amssymb}
\usepackage{mathrsfs}
\usepackage{empheq}
\usepackage{stmaryrd}
\usepackage{dsfont}
\usepackage{enumerate}
\usepackage{cancel}
\usepackage{comment}
\usepackage{pgfplots}
\usepackage{soul}
\usepackage{multirow}
\usetikzlibrary{patterns} 
\usepackage[applemac]{inputenc}
\usepackage{hyperref}
\usepackage{xcolor}

\definecolor{mygreen}{rgb}{0, 0.5, 0}

\newcommand{\1}{\mathds{1}}
\newcommand{\A}{\mathbb{A}}

\newcommand{\N}{\mathbb{N}}

\newcommand{\R}{\mathbb{R}}
\newcommand{\C}{\mathbb{C}}

\newcommand{\vsS}{\mathbb{S}}

\newcommand{\vH}{{\cal H}}

\newcommand{\vR}{{\cal R}}

\newcommand{\Ar}{\mathscr{A}}

\newcommand{\Dr}{\mathscr{D}}

\newcommand{\Lr}{\mathscr{L}}

\newcommand{\vphi}{\varphi}
\newcommand{\eps}{\varepsilon}
\newcommand{\dsp}{\displaystyle}
\newcommand{\ovl}{\overline}

\newcommand{\vlim}{\lim\limits}

\newcommand{\vmin}{\min\limits}

\newcommand{\vint}{\int\limits}

\newcommand{\inj}{\hookrightarrow}
\newcommand{\tends}{\longrightarrow}
\newcommand{\weak}{\rightharpoonup}
\newcommand{\wt}{\widetilde}

\newcommand{\loc}{\mathrm{loc}}

\newcommand{\rad}{\mathrm{rad}}

\newcommand{\co}{\mathrm{c}}
\renewcommand{\d}{\mathrm{d}}

\newcommand{\dist}{\mathrm{dist}}

\newcommand{\M}{\mathrm{max}}

\newcommand{\vi}{\mathrm{i}}
\newcommand{\w}{{\textsl w}}
\renewcommand{\le}{\leqslant}
\renewcommand{\ge}{\geqslant}
\renewcommand{\Re}{\mathrm{Re}}
\renewcommand{\Im}{\mathrm{Im}}

\newcommand{\p}{\prime}
\newcommand{\vect}{\overrightarrow}
\newcommand{\eqdef}{\stackrel{\mathrm{def}}{=}}
\DeclareMathOperator{\supp}{supp}

\numberwithin{equation}{section}

\newtheorem{thm}{Theorem}[section]
\newtheorem{prop}[thm]{Proposition}

\newtheorem{lem}[thm]{Lemma}
\newtheorem{sym}[thm]{Symmetry Property}

\theoremstyle{definition}
\newtheorem{rmk}[thm]{Remark}
\newtheorem{defi}[thm]{Definition}

\newtheorem{con}[thm]{Convention}

\newenvironment{proof*}{\noindent{\bf Proof.}}{\qed}
\newenvironment{vproof}[1]{\noindent{\bf Proof #1}}{\qed}

\title{\huge On the compactness of the support of solitary waves of the complex saturated nonlinear Schrödinger equation and related problems}
\author{\sc Pascal Bégout$^*$ and Jes\'us Ildefonso D{\'{\i}}az$^\dagger$}
\date{}

\begin{document}

\maketitle

\begin{center}
\textit{Dedicated to the memory of a master in nonlinear analysis: Ha\"{\i}m Brezis (1944-2024)}
\end{center}

\begin{gather*}
\begin{array}{cc}
	^*\text{Toulouse School of Economics}	&	\;^\dagger\text{Instituto de Matem\'atica Interdisciplinar}	\\
		\text{Université Toulouse Capitole}	&	\text{Universidad Complutense de Madrid}			\\
\text{Institut de Mathématiques de Toulouse}	&	\text{Plaza de las Ciencias, 3}						\\
		\text{1, Esplanade de l’Université}	&	\text{28040 Madrid, SPAIN}						\\
	\text{31080 Toulouse Cedex 6, FRANCE}	&
\bigskip \\
\text{
{\footnotesize E-mail\:: \href{mailto:Pascal.Begout@math.cnrs.fr}{\texttt{Pascal.Begout@math.cnrs.fr}}}}
&
\text{
{\footnotesize E-mail\:: \href{mailto:jidiaz@ucm.es}{\texttt{jidiaz@ucm.es}}}
}
\end{array}
\end{gather*}

\begin{abstract}
We study the vectorial stationary Schr\"{o}dinger equation $-\Delta u+a\,U+b\,u=F,$ with a saturated nonlinearity $U=u/|u|$ and with some complex coefficients $(a,b)\in\C^2.$ Besides the existence and uniqueness of solutions for the Dirichlet and Neumann problems, we prove the compactness of the support of the solution, under suitable conditions on $(a,b)$ and even when the source in the right hand side $F(x)$ is not vanishing for large values of $|x|.$ The proof of the compactness of the support uses a local energy method, given the impossibility of applying the maximum principle. We also consider the associate Schr\"{o}dinger-Poisson system when coupling with a simple magnetic field. Among other consequences, our results give a rigorous proof of the existence of  ``solitons with compact support" claimed, without any proof, by several previous authors.
\end{abstract}

{\let\thefootnote\relax\footnotetext{Pascal Bégout acknowledges funding from ANR under grant ANR-17-EUR-0010 (Investissements d’Avenir program)}}
{\let\thefootnote\relax\footnotetext{The research of J.~I.~D\'{\i}az was partially supported by the project PID2020-112517GB-I00 of the Spain State Research Agency (AEI)}}
{\let\thefootnote\relax\footnotetext{2020 Mathematics Subject Classification: 35J10 (35J88, 35Q51, 35R35, 35B30)}}
{\let\thefootnote\relax\footnotetext{Keywords: Schr\"{o}dinger equation, Schr\"{o}dinger-Poisson system, saturated nonlinear terms, solutions with compact support, local energy method, existence and uniqueness of solutions, solitons}}

\tableofcontents

\baselineskip .6cm

\section{Introduction}
\label{introduction}

We study the vectorial stationary Schr\"{o}dinger equation with a saturated nonlinearity 
\begin{gather}
\label{nls}
\begin{cases}
-\Delta u+a\,U+b\,u=F, \text{ in } H^{-1}(\Omega)+L^\infty(\Omega),	\medskip \\
U=\dfrac{u}{|u|}, \text{ almost everywhere in } \omega,
\end{cases}
\end{gather}
where $(a,b)\in\C^2,$ $\omega=\big\{x\in\Omega;u(x)\neq0\big\}$ and $\Omega$ is a subset of $\R^N$ whose boundary is $\Gamma,$ with homogeneous Dirichlet boundary condition
\begin{gather}
\label{dir}
 u_{|\Gamma}=0,
\end{gather}
or with homogeneous Neumann boundary condition
\begin{gather}
\label{neu}
 \dfrac{\partial u}{\partial\nu}_{|\Gamma}=0.
\end{gather}
We point out the important difference between the structure of the ``saturated nonlinearity" considered in this paper with respect to some related nonlinearities called by some authors as ``saturable nonlinearity" in which $U$ is replaced by some regularized approximations as, for instance, function $g_n$ below (see also the cases mentioned in the book by Agrawal  and Kivshar~\cite{ak}).

There are several different motivations to the study of such a coupled nonlinear system. The first one is associated to the consideration of non-Kerr law optical Schr\"{o}dinger equation with a saturated nonlinearity arising, for instance, in nonlinear optical media (see, e.g. Agrawal, G. P. and Kivshar~\cite{ak}, Biswas and Konar~\cite{MR2272971}) in which the saturated nonlinear term is understood in an approximated framework
\begin{gather*}
\vi\dfrac{\partial\psi}{\partial t}+\Delta\psi +a|\psi|^{-(1-m)}\psi=f(t,x),
\end{gather*}
(notice that our interest here corresponds to the choice $m=0)$ and we search for ``solitary wave solutions" of the form $\psi(t,x)=u(x)e^{\vi bt}$ (when $f(t,x)=e^{\vi bt}F(x)).$ This type of equations also arises in Quantum Mechanics and Hydrodynamics.

We mention that in a series of precedent papers (\cite{MR4340780,MR4725781}) we study the \textit{extinction in a finite time property} for the case of \textit{damped} non-Kerr law optical Schr\"{o}dinger equation with a saturated nonlinearity. Here our interest has a complementary nature since we are concentrated in the spatial behavior of solutions.

The different possible assumptions on the complex coefficients $(a,b)\in\C^2$ play a fundamental structural framework in the study of the system. In order to better advertise this fact, it is useful to rewrite the system in terms of the real components of solutions and data 
\begin{gather*}
\begin{cases}
u=u_R+\vi u_I,		&	F=F_R+\vi F_I,	\medskip \\ 
a=a_R+\vi a_I,		&	b=b_R+\vi b_I.
\end{cases}
\end{gather*}
Then, the sign of the components of coefficient $a$ is especially crucial for understanding the different nature of the coupled system
\begin{gather*}
\begin{cases}
-\Delta u_R+a_R\frac{u_R}{\sqrt{u_R^2+u_I^2}}-a_I\frac{u_I}{\sqrt{u_R^2+u_I^2}}+b_Ru_R-b_Iu_I=F_R,	\medskip \\ 
-\Delta u_I+a_R\frac{u_I}{\sqrt{u_R^2+u_I^2}}+a_I\frac{u_R}{\sqrt{u_R^2+u_I^2}}+b_Ru_I+b_Iu_R=F_I.
\end{cases}
\end{gather*}

Curiously enough, this kind of nonlinear coupled systems also arises (for different values of the coefficients) in the study of systems of Schr\"{o}dinger coupled equations with real coefficients: see, e.g. Ambrosetti~\cite{MR2447960} and Maia, Montefusco and Pellacci~\cite{MR3016511} for the case of non-linear terms of saturable type).

Systems of this type also arises in a very different setting: the study of compressed modes in Image Processing (see, e.g. Ozolins, Lai, Caflisch and Osher~\cite{olco}, Choukroun, Shtern, Bronstein and Kimmel~\cite{csbk} and its many references) when the system is understood as the corresponding Euler-Lagrange system associated to a variational functional involving the $L^{1}-$norm of the vectorial unknown. This last fact makes arise the presence of the saturated nonlinearity in the system of scalar equations.

Our study will deal also with the natural coupling of a Schr\"{o}dinger equation with a magnetic field in the spirit of many papers in the literature (see, Benci and Fortunato~\cite{MR1659454}, Colin and Watanabe~\cite{MR3639295} and their many references). As a matter of facts, we will consider merely the simpler case in which the coupling is reduced to the so-called Schr\"{o}dinger-Poisson system: looking for
\begin{gather}
\label{bcRN}
\begin{cases}
(u,U,\phi)\in H^1(\R^N)\cap L^1(\R^N)\times L^\infty(\R^N)\times\Dr^{1,2}(\R^N;\R),	\medskip \\
U \text{ is a saturated section associated to } u,								\medskip \\
\phi\,u\in L^2(\R^N) \text{ and } \phi\ge0 \text{ in  } \R^N,
\end{cases}
\end{gather}
solution to
\begin{gather}
\label{nlsSPRN}
\begin{cases}
-\Delta u+a\,U+b\,u+e\,\phi\,u=F, \text{ in } H^{-1}(\R^N)+L^\infty(\R^N),	\medskip \\
-\Delta\phi=\dfrac{e}2|u|^2, \text{ in } L^2(\R^N),
\end{cases}
\end{gather}
where the gauge potential $(\mathbf{A},\phi)$ is reduced to the case $\mathbf{A=0}$ and the real coefficient $e$ describes the strength of the iteration playing an important role in the study of the system. Note that if $e\le0$ then $\phi\le0$ and $e\phi\ge0,$ and if $e\ge 0$ then $\phi\ge 0$ and $e\phi\ge0$ (see Remark~\ref{rmkthmSPRN} below). So without loss of generality, we may assume that $e\ge0.$
\medskip

One of our main goals is to get sufficient conditions on the data implying that the corresponding ``solitary wave solution" $u$ has a compact support when $\Omega $ is unbounded. One of the motivations to consider this property comes from the fact that different authors claim the existence of solitary waves with compact support (called sometimes as ``solitons") but no rigorous proof of this fact was given until now (see, e.g., the many references collected in the monograph Biswas and Konar~\cite{MR2272971}). Here we will prove that this kind of ``solitons" exists in the presence of saturated non-linear terms, under suitable conditions on the source term $F.$

This kind of results on solutions with compact support was initiated with the pioneering paper by H.~Brezis~\cite{MR0481460} who constructed suitable super and subsolutions with compact support for some unilateral elliptic problems: the a priori unknown boundary of the support of the solution is a ``free boundary" of the problem. After that, the result was extended to many other equations for which the maximum principle holds (see, e.g. the exposition made in \cite{MR853732}). The study of this property for problems in which the maximum principle may fail was the main object of the monograph by Antontsev, D\'iaz and Shmarev~\cite{MR2002i:35001} where different local energy methods were proposed to get such type of qualitative properties. This approach was improved and applied to sublinear Schr\"{o}dinger equations in \cite{MR2876246} and \cite{MR3190983}. The adaptation of the energy methods to the case of parabolic problems involving some possible multivalued nonlinear term was carried out in \cite{MR853732}. The approach we will follow in this paper is different and simpler than the method used in \cite{MR2466410}. One of the unexpected facts we will prove, in comparison with the results obtained by the authors for suitable sublinear Schr\"{o}dinger equations, is that the compact support may hold even if the right hand side term $F$ is never zero in the whole domain $\Omega.$ As mentioned before, this case was considered by the authors when dealing with a completely different property: the extinction of the solution in a finite time for the evolution damped saturated Schr\"{o}dinger equations, nevertheless the required assumptions on the coefficients $(a,b)\in\C^2$ are not the same. Among many general cases, we will prove that the compact support property can be satisfied in the case in which the coefficients of the real or the imaginary parts of the unknown are not both of ``pure absorption type" (case which corresponds to the assumptions $a_R>0,$ $b_R>0$ and $a_I=b_I=0):$ see assumption \eqref{ab} below. This also explain why the uniqueness of solutions is a delicate subject requiring some suitable assumptions on $(a,b)\in\C^2.$

A very particular statement, dealing with the associated evolution Schr\"{o}dinger nonlinear problem, can be presented here as an application of Theorems~\ref{thmexi}, \ref{thmuni}, and \ref{thmsolcomRN}, and \cite[Proposition2.6]{MR4725781} below.
\begin{thm}
\label{thmA}
Let $\lambda\in\{1,\vi,-\vi\},$ $\mu>0$ and $F\in H^{-1}(\R^N).$ Assume that $\|F\|_{L^\infty(K^\co)}$ is small enough for some compact subset $K$ of $\R^N.$ Let $u_0\in H^1(\R^N)$ be the unique global weak solution  $($see Definition~$\ref{defsol}$ below$)$ to
\begin{gather*}
-\Delta u_0+\mu u_0+\lambda\dfrac{u_0}{|u_0|}=-F, \text{ in } H^{-1}(\R^N)+L^\infty(\R^N).
\end{gather*}
For $(t,x)\in\R\times\R^N,$ set $u(t,x)=u_0(x)e^{\vi\mu t}.$ Then, $u\in C^\infty(\R;H^1(\R^N))$ is a solution to
\begin{gather*}
\begin{cases}
\vi\dfrac{\partial u}{\partial t}+\Delta u=\lambda\dfrac{u}{|u|}+F(x)e^{\vi\mu t},	&	\text{in } \R\times\R^N,	\medskip \\
u(0)=u_0, 														&	\text{in } \R^N.
\end{cases}
\end{gather*}
Furthermore, $\supp u(t)=\supp u_0$ is compact, for any $t\in\R.$ Finally, if $\lambda=-\vi$ and $F\in L^2(\R^N)$ then the solution $u$ is unique.
\end{thm}
The organization of this paper is as follows: below, we present a short collection of the notations used in the paper. The statements on the
existence and uniqueness of solutions is given in Section~\ref{seceu} under suitable assumptions on the complex coefficients $(a,b)\in\C^2$ which are
geometrically represented in Section~\ref{geoint} for the reader convenience. The statemets of the results on the compactness of the support of the solution is offered in Section~\ref{secsolcom}, after a previous section (Section~\ref{secine}) devoted to the presentation of some crucial inequalities which will be used in the energy method used in this paper. Section~\ref{secproofeu} contains the proofs of the existence and uniqueness of solutions while the proofs of the spatial localization inequalities and the statements on solutions compactly supported are given in Sections~\ref{secproofine} and~\ref{seproofsolcom}, respectively. Finally, the Schr\"{o}dinger-Poisson system is considered in Section~\ref{SP} where the compactness of the support of the component $u$ is proved under suitable conditions on the data.

For a complex number $z,$ we denote by $\ovl z,$ $\Re(z)$ and $\Im(z),$ its conjugate, real and imaginary part, respectively, and $\vi^2=-1.$ For $p\in[1,\infty],$ $p^\prime$ is the conjugate of $p$ defined by $\frac{1}{p}+\frac{1}{p^\prime}=1.$ Let $\Omega$ be an open subset of $\R^N.$ Unless specified, all functions are complex-valued $(H^1(\Omega)\eqdef H^1(\Omega;\C),$ etc) and all the vector spaces are considered over the field $\R.$ For a (real) Banach space $X,$ we denote by $X^\star\eqdef\Lr(X;\R)$ its topological dual and by $\langle\: . \; , \: . \:\rangle_{X^\star,X}\in\R$ the $X^\star-X$ duality product. When $X$ (respectively, $X^\star)$ is endowed of the weak topology $\sigma(X,X^\star)$ (respectively, the weak$\star$ topology $\sigma(X^\star,X)),$ it is denoted by $X_\w$ (respectively, by $X_{\w\star}).$ Auxiliary positive constants will be denoted by $C$ and may change from a line to another one. Also for positive parameters $a_1,\ldots,a_n,$ we shall write $C(a_1,\ldots,a_n)$ to indicate that the constant $C$ depends only and continuously on $a_1,\ldots,a_n.$ If $X$ and $E$ are two Banach spaces then $X\inj E$ means that $X\subset E$ and that the identity function $i:X\tends E$ is (injective and) continuous. If $A$ is a subset of $\R^N$ then $A^\co$ denotes its complement, and $A\setminus B=A\cap B^\co.$
\medskip \\
Let us recall some well-known results that we shall often use without referring to them. The first one may be proved with Egorov' Theorem (see also Strauss~\cite{MR306715}). Let $\Omega$ be an open subset of $\R^N,$ $u\in\C^\Omega,$ $1<p\le\infty,$ and let $(u_n)_{n\in\N}\subset L^p(\Omega)$ be bounded. If $u_n\xrightarrow[n\to\infty]{\text{a.e.\;in }\Omega}u$ then $u\in L^p(\Omega)$ and $u_n\xrightarrow[n\to\infty]{L^q(\Omega^\p)}u,$ for any $1\le q<p$ and subset $\Omega^\p\subseteq\Omega$ of finite measure. In addition, $u_n\underset{n\to\infty}{\overset{L^p(\Omega)_\w}{-\!\!\!-\!\!\!-\!\!\!-\!\!\!-\!\!\!\weak}}u,$ if $p<\infty,$ and $u_n\underset{n\to\infty}{\overset{L^\infty(\Omega)_{\w\star}}{-\!\!\!-\!\!\!-\!\!\!-\!\!\!-\!\!\!-\!\!\!-\!\!\!\weak}}u,$ if $p=\infty.$ The second one comes from results for the weak topology, the Rellich-Kondrachov Theorem and a classical result of integration (H.~Brezis~\cite[Theorem~4.9, p.94]{MR2759829}). Let $(u_n)_{n\in\N}\subset H^1_0(\Omega)$ be bounded. Then there exists $u\in H^1_0(\Omega)$ such that, up to a subsequence, $u_n\underset{n\to\infty}{\overset{H^1_0(\Omega)_\w}{-\!\!\!-\!\!\!-\!\!\!-\!\!\!-\!\!\!-\!\!\!\weak}}u,$ $u_n\xrightarrow[n\to\infty]{L^2_\loc(\Omega)}u,$ and $u_n\xrightarrow[n\to\infty]{\text{a.e.\;in }\Omega}u.$ There is a similar statement in $H^1(\Omega).$ The third result is the following. If $X$ and $E$ are two Banach spaces such that $X\inj E$ with dense embedding then $E^\star\inj X^\star,$ and
\begin{gather*}
\forall F\in E^\star, \; \forall u\in X, \; \langle F,u\rangle_{X^\star,X}=\langle F,u\rangle_{E^\star,E}.
\end{gather*}
By the Riesz representation Theorem, we have for any $p\in[1,\infty),$
\begin{gather*}
\forall F\in L^{p^\p}(\Omega), \; \forall u\in L^p(\Omega), \; \langle F,u\rangle_{L^{p^\p}(\Omega),L^p(\Omega)}=\Re\vint_\Omega F(x)\ovl{u(x)}\d x.
\end{gather*}
In particular, this implies that we shall always identify $L^2(\Omega)$ with its topological dual. Finally, if $A_1$ and $A_2$ are two Banach spaces such that $A_1,A_2\subset\vH$ for some Hausdorff topological vector space $\vH,$ and if $A_1\cap A_2$ is dense in both $A_1$ and $A_2$ then $A_1\cap A_2$ and $A_1+A_2$ are Banach spaces, whose norms are
\begin{gather*}
\|a\|_{A_1\cap A_2}=\max\big\{\|a\|_{A_1},\|a\|_{A_2}\big\}	\; \text{ and } \;
\|a\|_{A_1+A_2}=\inf_{\left\{\substack{a=a_1+a_2 \hfill \\ (a_1,a_2)\in A_1\times A_2}\right.}\Big(\|a_1\|_{A_1}+\|a_2\|_{A_2}\Big),
\end{gather*}
respectively, and $\big(A_1\cap A_2\big)^\star=A_1^\star+A_2^\star.$ This justifies the identity~\eqref{defsol21} below. For more details, see Trèves~\cite{MR2296978}, Bergh and Löfström~\cite{MR0482275}, and~\cite{MR4521439}.
\medskip \\
We end our reminders by the space $\Dr^{1,2}(\R^N),$ with $N\ge3.$ Setting $|u|_2=\|\nabla u\|_{L^2(\R^N)},$ then $\Dr^{1,2}(\R^N)$ is the completion of $\Dr(\R^N)$ with respect to the norm $|\:.\:|_2.$ We have that $\Dr^{1,2}(\R^N)$ is a reflexive separable Banach space, the embedding $\Dr^{1,2}(\R^N)\inj L^2_\loc(\R^N)$ is compact, and
\begin{gather*}
\Dr^{1,2}(\R^N)=\big\{u\in L^\frac{2N}{N-2}(\R^N);\nabla u\in L^2(\R^N;\C^N)\big\}.
\end{gather*}
Finally, there exists $C=C(N)$ such that for any $u\in\Dr^{1,2}(\R^N;\R),$ $\|u\|_{L^\frac{2N}{N-2}(\R^N)}\le C\|\nabla u\|_{L^2(\R^N)}.$ For more details on these spaces, see for instance Galdi~\cite[Chapter~II]{MR2808162} (Sections~II.6--II.10), and Lieb and Loss~\cite{MR1817225} (Theorem~8.6, p.208, and Corollary~8.7, p.212).

\section{Existence and uniqueness}
\label{seceu}

\begin{defi}
\label{defsol}
Let $\Omega\subseteq\R^N$ be an open subset, $(a,b)\in\C^2,$ and
\begin{gather}
\label{phi}
\phi\in L^\infty(\Omega)+L^{p_\phi}(\Omega), \text{ where}		\\
\label{pphi}
p_\phi=
\begin{cases}
2,								&	\text{if } N=1,	\\
2+\kappa, \text{ for some } \kappa>0,	&	\text{if } N=2,	\\
N,								&	\text{if } N\ge3.
\end{cases}
\end{gather}
Although not necessary in this definition, it will sometimes be convenient to assume that $\phi$ verifies
\begin{gather}
\label{phip}
\phi:\Omega\tends\R \; \text{ and } \; \phi\ge0, \text{ a.e.\,in } \Omega.
\end{gather}
\begin{enumerate}
\item
\label{defsol1}
Let $u\in L^1_\loc(\Omega).$ A function $U$ is said to be a \textit{saturated section} associated to $u$ if it satisfies
\begin{gather}
\label{defsol11}
U\in L^\infty(\Omega) \text{ and } \|U\|_{L^\infty(\Omega)}\le1,	\\
\label{defsol12}
U=\frac{u}{|u|}, \text{ almost everywhere in } \omega,
\end{gather}
where $\omega=\big\{x\in\Omega;u(x)\neq0\big\}.$
\item
\label{defsol2}
Let $F\in H^{-1}(\Omega)+L^\infty(\Omega).$ We shall say that a function $u$ is a \textit{global weak solution} to
\begin{gather}
\label{nlsg}
-\Delta u+a\,U+b\,u+\phi\,u=F, \text{ in } H^{-1}(\Omega)+L^\infty(\Omega),
\end{gather}
with boundary condition~\eqref{dir}, if $u\in H^1_0(\Omega)\cap L^1(\Omega),$ $U$ is a saturated section associated to $u,$ and if
\begin{multline}
\label{defsol21}
\langle\nabla u,\nabla v\rangle_{L^2(\Omega),L^2(\Omega)}+\langle a\,U,v\rangle_{L^\infty(\Omega),L^1(\Omega)}
												+\langle b\,u,v\rangle_{L^2(\Omega),L^2(\Omega)}	\\
+\langle\phi\,u,v\rangle_{L^2(\Omega),L^2(\Omega)}=\langle F,v\rangle_{X^\star,X},
\end{multline}
for any $v\in H^1_0(\Omega)\cap L^1(\Omega),$ where $X=H^1_0(\Omega)\cap L^1(\Omega).$
\item
\label{defsol3}
Assume that $\Omega$ has a finite measure and a Lipschitz continuous boundary. Let $F\in H^1(\Omega)^\star.$ We shall say that a function $u$ is a \textit{global weak solution} to \eqref{nlsg} with boundary condition \eqref{neu} if $u\in H^1(\Omega),$ $U$ is a saturated section associated to $u,$ and if $(u,U,\phi)$ satisfies \eqref{defsol21} for any $v\in H^1(\Omega),$ where $X=H^1(\Omega).$
\end{enumerate}
Sometimes, we shall write $(u,U)$ or $(u,U,\phi)$ to designate a solution with the obvious meanings.
\end{defi}

\begin{rmk}
\label{rmkphiu}
Let $\phi=\phi_1+\phi_2\in L^\infty(\Omega)+L^{p_\phi}(\Omega),$ where $p_\phi$ is given by \eqref{pphi}. If $u\in H^1_0(\Omega)$ (or if $u\in H^1(\Omega)$ and $\Omega$ has a Lipschitz continuous boundary) then it follows from Hölder's inequality and Sobolev' embedding that,
\begin{gather*}
\|\phi u\|_{L^2(\Omega)}\le\|\phi_1\|_{L^\infty(\Omega)}\|u\|_{L^2(\Omega)}+C\|\phi_2\|_{L^{p_\phi}(\Omega)}\|u\|_{H^1_0(\Omega)},
\end{gather*}
from which we deduce,
\begin{gather}
\label{rmkphiu1}
\|\phi u\|_{L^2(\Omega)}\le C\|\phi\|_{L^\infty(\Omega)+L^{p_\phi}(\Omega)}\|u\|_{H^1_0(\Omega)},
\end{gather}
where $C=C(N,\kappa).$ It follows that the $L^2-L^2$ duality product in \eqref{defsol21} involving $\phi u$ makes sense, and we have
\begin{gather}
\label{rmkphiu2}
|\langle\phi\,u,v\rangle_{L^2(\Omega),L^2(\Omega)}|\le C\|\phi\|_{L^\infty(\Omega)+L^{p_\phi}(\Omega)}\|u\|_{H^1_0(\Omega)}\|v\|_{L^2(\Omega)},
\end{gather}
for any $v\in L^2(\Omega),$ where $C=C(N,\kappa).$
\end{rmk}

\begin{con}
\label{conneu}
Throughout this paper, $\Omega$ denotes any open subset of $\R^N,$ and $(a,b)$ is a pair of complex numbers. When a function will be said to satisfy the boundary condition \eqref{neu}, it will always be assumed that $\Omega$ has a finite measure and a Lipschitz continuous boundary.
\end{con}

\noindent
Let,
\begin{gather}
\label{A}
\A=\C\setminus\big\{z\in\C; \Re(z)\le0 \text{ and } \Im(z)=0\big\},	
\end{gather}
\begin{gather}
\label{ab}
(a,b)\in\A\times\A \quad \text{ and } \quad
\begin{cases}
\Im(a)\Im(b)\ge0,										\medskip \\
\text{ or }												\medskip \\
\Im(a)\Im(b)<0 \; \text{ and } \; \Re(b)>\dfrac{\Im(b)}{\Im(a)}\Re(a).
\end{cases}
\end{gather}

\begin{thm}[\textbf{Null solution}]
\label{thmunull}
Let $(a,b)$ satisfy~\eqref{ab}, with eventually $b=0,$ let $\phi$ satisfy~\eqref{phi}--\eqref{phip}, and let $F\in L^\infty(\Omega)$ with $\|F\|_{L^\infty(\Omega)}\le|a|.$ Then there exists $M=M(|a|,|b|)$ such that if $\|F\|_{L^\infty(\Omega)}\le\frac1{M}$ the unique global weak solution $(u,U)$ to~\eqref{nlsg} with boundary condition~\eqref{dir} or \eqref{neu} is given by,
\begin{gather}
\label{thmunull1}
u=0 \; \text{ and } \; U=\frac1aF,
\end{gather}
almost everywhere in $\Omega.$
\end{thm}

\begin{rmk}
\label{rmkthmunull}
Let $(u,U)$ be defined by \eqref{thmunull1}. It is easy to check that if $\|F\|_{L^\infty(\Omega)}\le|a|$ then $(u,U)$ is always a solution. If $\|F\|_{L^\infty(\Omega)}>|a|$ or if $F\notin L^\infty(\Omega)$ then $(u,U)$ is never a solution, in the sense of Definition~\ref{defsol}. Indeed, $U$ does not satisfies \eqref{defsol11}. This is trivially obtained by the identity $\|U\|_{L^\infty(\Omega)}=\frac1{|a|}\|F\|_{L^\infty(\Omega)}.$
\end{rmk}

\begin{thm}[\textbf{Existence}]
\label{thmexi}
Let $(a,b)$ satisfy~\eqref{ab}, and let $\phi$ satisfy~\eqref{phi}--\eqref{phip}. Then, for any $F\in H^{-1}(\Omega)$ $($respectively, $F\in H^1(\Omega)^\star),$ there exists at least one global weak solution to~\eqref{nlsg} with boundary condition~\eqref{dir} $($respectively,~\eqref{neu}$).$ In addition, Symmetry Property~$\ref{sym}$ below holds.
\end{thm}

\begin{sym}
\label{sym}
If, furthermore, there exists $\vR\in SO_N(\R)$ such that for almost every $x\in\Omega,$ $\vR x\in\Omega,$ $F(\vR x)=F(x)$ and $\phi(\vR x)=\phi(x)$ then we may construct a solution $u$ which also satisfies $u(\vR x)=u(x),$ for almost every $x\in\Omega.$ When $N=1,$ if $\Omega$ is symmetric with respect to the origin and if $F$ and $\phi$ are odd functions then $u$ is also an odd function.
\end{sym}

\noindent
Here and in what follows, $SO_N(\R)$ denotes the special orthogonal group of $\R^N.$ 

\begin{thm}[\textbf{Uniqueness}]
\label{thmuni}
Let $\phi$ satisfy~\eqref{phi}--\eqref{pphi}, and let $(a,b)\in\C^2$ be such that $a\neq0,$ $\Re(a)\ge0$ and $\Re(a\ovl b)+\Re(a\ovl\phi)\ge0,$ a.e.\,in $\Omega.$ If $\Re(a\ovl b)+\Re(a\ovl\phi)=0,$ on a set of positive measure, then assume further that one of the three following conditions is satisfied.
\begin{enumerate}
\item
$\Re(b)+\Re(\phi)>0,$ a.e.\,in $\Omega.$
\item
$\Im(b)+\Im(\phi)>0,$ a.e.\,in $\Omega.$
\item
$\Im(b)+\Im(\phi)<0,$ a.e.\,in $\Omega.$
\end{enumerate}
If $F\in H^{-1}(\Omega)+L^\infty(\Omega)$ $($respectively, $F\in H^1(\Omega)^\star)$ and if $(u_1,U_1)$ and $(u_2,U_2)$ are two global weak solutions to~\eqref{nlsg} with boundary condition~\eqref{dir} $($respectively,~\eqref{neu}$)$ then $u_1=u_2$ and $U_1=U_2.$
\end{thm}

\begin{rmk}
\label{rmkthmuni}
Assume that $\phi=0.$ Then the assumptions simply become $a\neq0,$ $\Re(a)\ge0,$ $\Re(a\ovl b)\ge0$ and $b\in\A.$ Note that in any case, $a\in\A.$ In~\cite{MR3315701}, the authors study, among others, $-\Delta u+a|u|^{-(1-m)}u+bu=F,$ with $0<m<1.$ Uniqueness is ensured if $a\neq0,$ $\Re(a)\ge0$ and $\Re(a\ovl b)\ge0$ (\cite[Theorem~2.10]{MR3315701}). Here, the equation we study corresponds to $m=0,$ and we have to add an assumption to $b,$ namely $b\in\A.$ This is due to the fact that the nonlinearity $|u|^{-(1-m)}u$ becomes $\frac{u}{|u|},$ whose the behavior near $0$ is unknown.
\end{rmk}

\begin{thm}[\textbf{A priori bound}]
\label{thmbound}
Let $(a,b)$ satisfy~\eqref{ab}, let $\phi$ satisfy~\eqref{phi}--\eqref{phip}, and let $F\in H^{-1}(\Omega)$ $($respectively, $F\in H^1(\Omega)^\star).$ Then there exists $M=M(|a|,|b|)$ such that any global weak solution $u$ to~\eqref{nlsg} with boundary condition~\eqref{dir} $($respectively,~\eqref{neu}$)$ satisfies
\begin{gather}
\label{thmbound11}
\|u\|_X^2+\|u\|_{L^1(\Omega)}+\int_\Omega\phi|u|^2\d x\le M\|F\|_{X^\star}^2,
\end{gather}
where $X=H^1_0(\Omega)$ $($respectively, $X=H^1(\Omega)).$
\end{thm}

\section{Geometric interpretation of existence and uniqueness of complex parameters}
\label{geoint}

\textbf{Existence.}
Throughout this paper, the existence of a solution of the equations we study is ensured if $(a,b)$ satisfies \eqref{ab} (Theorem~\ref{thmexi}). Seeing the complex numbers $a$ and $b$ as points in the Euclidean plane, the geometric representation of this condition is the following: $[a,b]\cap\Ar=\emptyset,$ where $\Ar$ is the geometric representation of $\A^\co:$ the half red line where $\Re(z)\le0$ and $\Im(z)=0.$ See Figures~1 and 2 below.
\medskip

\begin{tikzpicture}
    \fill[fill=gray!20] (0,0) -- (3,-3) -- (3,3) -- (-3,3) -- (-3,0) -- cycle;
    
    \fill[fill=red!20] (0,0) -- (3,-3) -- (-3,-3) -- (-3,0) -- cycle;
    
    \draw[red] (-3,0) -- (0,0);
    \draw[->] (0,0) -- (3,0) node[right] {$\Re(z)$};
    \draw[->] (0,-3) -- (0,3) node[above] {$\Im(z)$};
    
    \node[circle,color=red,fill,inner sep=1.5pt,label={below left:$0$}] at (0,0) {};
    
    \node[circle,fill,inner sep=1.5pt,label={above right:$a$}] at (-1.5,1.5) {};
    
    \draw[dashed] (-1.5,1.5) -- (0,0);
    
    \draw[red] (0,0) -- (3,-3);
    
    \node[below right] at (1.3,-1) {$\scriptstyle\Re(z)=\frac{\Re(a)}{\Im(a)}\Im(z)$};
            
    \filldraw[gray!20] (-3,-4) rectangle (-2,-3.5) node[midway, right, black, xshift=15pt] {\footnotesize Admissible values for $b$ with respect to $a$};
        
    \filldraw[red!20] (-3,-4.75) rectangle (-2,-4.25) node[midway, right, black, xshift=15pt] {\footnotesize Forbidden values for $b$ with respect to $a$};

    \node[below] at (0,-5) {\small Figure 1: Existence with $\Im(a)\neq0$};
\end{tikzpicture}
\hspace{1cm}
\begin{tikzpicture}
    \fill[fill=red!20] (-3,0) rectangle (0,3);
    
    \fill[fill=gray!20] (0,0) -- (-3,0) -- (-3,-3) -- (3,-3) -- (3,3) -- (0,3) -- cycle;
    
    \draw[red] (-3,0) -- (0,0);
    \draw[->] (0,0) -- (3,0) node[right] {$\Re(z)$};
    \draw[] (0,-3) -- (0,0);
    \draw[red][->] (0,0) -- (0,3) node[above] {\color{black}$\Im(z)$};
    
    \node[circle,color=red,fill,inner sep=1.5pt,label={below left:$0$}] at (0,0) {};
    
    \node[circle,fill,inner sep=1.5pt,label={above right:$a$}] at (0,-2) {};
    
    \filldraw[gray!20] (-3,-4) rectangle (-2,-3.5) node[midway, right, black, xshift=15pt] {\footnotesize Admissible values for $b$ with respect to $a$};
        
    \filldraw[red!20] (-3,-4.75) rectangle (-2,-4.25) node[midway, right, black, xshift=15pt] {\footnotesize Forbidden values for $b$ with respect to $a$};

    \node[below] at (0,-5) {\small Figure 2: Existence with $\Re(a)=0$};
\end{tikzpicture}

\noindent
\textbf{Uniqueness.}
Assume $ab\neq0$ and $\phi=0.$ Uniqueness is ensured if $\Re(a)\ge0,$ $\Re(a\ovl b)\ge0$ and $b\in\A$ (Theorem~\ref{thmuni} and Remark~\ref{rmkthmuni}). Now, let us see a complex number $z$ as a vector of the Euclidean plane:
$
\vect z=
 \left(
  \begin{array}{c}
   \Re(z) \\
   \Im(z)
  \end{array}
 \right).
$
It follows that $\Re\left(a\ovl b\right)=\vect a.\vect b\ge0$ (here, $.$ denotes the scalar product between two vectors) means that $\left|\angle(\vect a,\vect b)\right|\le\dfrac{\pi}{2}\rad$ (see Figures~3 and 4 below).

\begin{tikzpicture}
    \fill[fill=gray!20] (-3,1.5) -- (3,-1.5) -- (3,3) -- (-3,3) -- cycle;
    
    \fill[fill=red!20] (-3,1.5) -- (3,-1.5) -- (3,-3) -- (-3,-3) -- cycle;
    
    \draw[->] (-3,0) -- (3,0) node[right] {$\Re(z)$};
    \draw[->] (0,-3) -- (0,3) node[above] {$\Im(z)$};
    
    \node[circle,color=red,fill,inner sep=1.5pt,label={below left:$0$}] at (0,0) {};
    
    \node[circle,fill,inner sep=1.5pt,label={above right:$a$}] at (1,2) {};
    
    \draw[dashed] (1,2) -- (0,0);
    
    \draw[black] (-3,1.5) -- (3,-1.5);
    
    \draw[black] (0.2,-0.1) -- (0.3,0.1);
    \draw[black] (0.1,0.2) -- (0.3,0.1);
    
    \node[below right] at (-3,1.85) {$\scriptstyle\Re(a\ovl z)=0$};
            
    \filldraw[gray!20] (-3,-4) rectangle (-2,-3.5) node[midway, right, black, xshift=15pt] {\footnotesize Admissible values for $b$ with respect to $a$};
        
    \filldraw[red!20] (-3,-4.75) rectangle (-2,-4.25) node[midway, right, black, xshift=15pt] {\footnotesize Forbidden values for $b$ with respect to $a$};

    \node[below] at (0,-5) {\small Figure 3: Uniqueness with $\Re(a)>0$};
\end{tikzpicture}
\hspace{1cm}
\begin{tikzpicture}
    \fill[fill=red!20] (-3,0) rectangle (3,3);
    
    \fill[fill=gray!20] (-3,-3) rectangle (3,0);
    
    \draw[red] (-3,0) -- (0,0);
    \draw[->] (0,0) -- (3,0) node[right] {$\Re(z)$};
    \draw[->] (0,-3) -- (0,3) node[above] {$\Im(z)$};
    
    \node[circle,color=red,fill,inner sep=1.5pt,label={below left:$0$}] at (0,0) {};
    
    \node[circle,fill,inner sep=1.5pt,label={above right:$a$}] at (0,-2) {};
    
    \filldraw[gray!20] (-3,-4) rectangle (-2,-3.5) node[midway, right, black, xshift=15pt] {\footnotesize Admissible values for $b$ with respect to $a$};
        
    \filldraw[red!20] (-3,-4.75) rectangle (-2,-4.25) node[midway, right, black, xshift=15pt] {\footnotesize Forbidden values for $b$ with respect to $a$};

    \node[below] at (0,-5) {\small Figure 4: Uniqueness with $\Re(a)=0$};
\end{tikzpicture}

\section{Inequalities for spatial localization}
\label{secine}

In this section and in what follows, $r_+=\max\{0,r\}$ denotes the positive part of the real number $r.$ For $x_0\in\R^N$ and $r>0,$ $B(x_0,r)$ is the open ball of $\R^N$ of center $x_0$ and radius $r,$ $\vsS(x_0,r)$ is its boundary, and $\ovl B(x_0,r)$ is its closure. Finally, $\sigma$ is the surface measure on a sphere.

\begin{thm}
\label{thmsta}
Let $N\in\N,$ $x_0\in\R^N,$ $\rho_0>0$ and $u\in H^1\big(B(x_0,\rho_0)\big).$ There exists $C=C(N)$ such that, if $u$ satisfies the inequality
\begin{gather}
\label{thmsta1}
\|\nabla u\|_{L^2(B(x_0,\rho))}^2+\|u\|_{L^1(B(x_0,\rho))}\le M\left|\dsp\int_{\vsS(x_0,\rho)}u\ovl{\nabla u}.\frac{x-x_0}{|x-x_0|}\d\sigma\right|,
\end{gather}
for almost every $\rho\in(0,\rho_0)$ and for some $M>0,$ then $u=0,$ almost everywhere in $B(x_0,\rho_\M),$ where
\begin{gather}
\label{thmsta2}
\rho_\M^{N+2}=\left(\rho_0^{N+2}-CM^2\max\left\{\rho_0^{N+1},1\right\}\vmin_{\tau\in\left(\frac12,1\right]}\left\{\frac{E(\rho_0)^{\gamma(\tau)}
\max\{b(\rho_0)^{\mu(\tau)},b(\rho_0)^{1-\gamma(\tau)}\}}{2\tau-1}\right\}\right)_+,
\end{gather}
and where,
\begin{gather*}
\begin{array}{lll}
E(\rho_0)=\|\nabla u \|_{L^2(B(x_0,\rho_0))}^2,		&	&	b(\rho_0)=\|u\|_{L^1(B(x_0,\rho_0))},	\\
\gamma(\tau)=\dfrac{2\tau-1}{N+2},				&	&	\mu(\tau)=\dfrac{2(1-\tau)}{N+2},
\end{array}
\end{gather*}
for any $\tau\in\left(\frac12,1\right].$
\end{thm}

\noindent
Below, we give a sufficient condition which ensures that $\rho_\M=\rho_0.$

\begin{thm}
\label{thmstaF}
Let $M>0,$ $L>0,$ $\rho_1>\rho_0>0,$ $F\in L^2\big(B(x_0,\rho_1)\big)$ and $u\in H^1\big(B(x_0,\rho_1)\big).$ Assume that for almost every $\rho\in(0,\rho_1),$
\begin{gather}
\label{thmstaF1}
\|u\|_{H^1(B(x_0,\rho))}^2+\|u\|_{L^1(B(x_0,\rho))}
\le M\left(\left|\int_{\vsS(x_0,\rho)}u\ovl{\nabla u}.\frac{x-x_0}{|x-x_0|}\d\sigma\right|+\int_{B(x_0,\rho)}|F(x)u(x)|\d x\right).
\end{gather}
Then, there exist $E_\star=E_\star(M,L,N,\rho_1,\rho_0)$ and $\eps_\star=\eps_\star(M,L,N,\rho_0,\rho_1)$ such that if $\|u\|_{L^1(B(x_0,\rho_1))}\le L,$ $\|\nabla u\|_{L^2(B(x_0,\rho_1))}^2\le E_\star$ and if
\begin{gather}
\label{thmstaF2}
\|F\|_{L^2(B(x_0,\rho))}^2\le\eps_\star\big((\rho-\rho_0)_+\big)^{N+2},
\end{gather}
for any $\rho\in(0,\rho_1),$ then $u=0,$ almost everywhere in $B(x_0,\rho_0).$
\end{thm}

\section{Solutions compactly supported}
\label{secsolcom}

Throughout this section, we assume that $(a,b)$ and $\phi$ satisfy \eqref{ab} and \eqref{phi}--\eqref{phip}, respectively.

\begin{thm}
\label{thmsolcomRN}
Let $F\in H^{-1}(\R^N),$ and let $u\in H^1(\R^N)\cap L^1(\R^N)$ be any global weak solution to~\eqref{nlsg}. There exists $M=M(|a|,|b|)$ such that if there exists a compact subset $K$ of $\R^N$ such that $F_{|K^\co}\in L^\infty(K^\co)$ and $\|F\|_{L^\infty(K^\co)}\le\frac1M,$ then $\supp u$ is compact.
\end{thm}

\begin{thm}
\label{thmsolcom}
Let $K\subset\Omega$ be any compact subset of $\R^N.$ Let $F\in H^{-1}(\Omega)$ be such that $F_{|\Omega\setminus K}\in L^\infty(\Omega\setminus K)$ and let $u$ be any global weak solution to~\eqref{nlsg} with boundary condition \eqref{dir}. Then there exist $M=M(|a|,|b|)$ and $\eps_\star=\eps_\star(\dist(K,\Gamma))$ such that, for any $\eps\in(0,\eps_\star),$ there exists $\delta=\delta(\eps,|a|,|b|,N)$ satisfying that if $\|F\|_{H^{-1}(\Omega)}<\delta$ and $\|F\|_{L^\infty(\Omega\setminus K)}\le\frac1M$ then $\supp u\subset K(\eps)\subset\Omega,$ where
\begin{gather}
\label{K}
K(\eps)=\Big\{x\in\R^N;\dist(x,K)\le\eps\Big\}.
\end{gather}
A similar statement holds for the boundary condition~\eqref{neu}.
\end{thm}

\begin{rmk}
\label{rmkthmsolcom}
Here are some comments about Theorems~\ref{thmsolcomRN}--\ref{thmsolcom}.
\begin{enumerate}
\item
\label{rmkthmsolcom1}
The solutions involved exist with help of Theorem~\ref{thmexi}.
\item
\label{rmkthmsolcom2}
By Theorem~\ref{thmunull}, if $\|F\|_{L^\infty(\Omega)}$ is small enough then the solution is $u=0,$ in $\Omega,$ and the above theorems do not give any new information. On the other hand, if $F_{|K}\notin L^\infty(K)$ or if $\|F\|_{L^\infty(K)}$ is large enough then $u=0,$ over $\Omega,$ can not be a solution (Remark~\ref{rmkthmunull}). It follows that the above results provide an additional qualitative property about the solutions.
\end{enumerate}
\end{rmk}

\section{Existence and uniqueness: the proofs}
\label{secproofeu}

Let $(a,b)$ satisfy~\eqref{ab}, with eventually $b=0,$ and let $\phi$ satisfy~\eqref{phi}--\eqref{phip}. Let $u$ be a global weak solution to \eqref{nlsg} with boundary condition \eqref{dir} or \eqref{neu}, with $F\in Y^\star.$ Here, $Y=H^1_0(\Omega)\cap L^1(\Omega),$ if $u$ satisfies \eqref{dir}, and $Y=H^1(\Omega),$ if $u$ satisfies \eqref{neu}. Choosing $u$ and $\vi u$ as test functions, we obtain
\begin{gather}
\label{predem1}
\|\nabla u\|_{L^2(\Omega)}^2+\Re(a)\|u\|_{L^1(\Omega)}+\Re(b)\|u\|_{L^2(\Omega)}^2+\int_\Omega\phi|u|^2\d x=\langle F,u\rangle_{Y^\star,Y},	\\
\label{predem2}
\Im(a)\|u\|_{L^1(\Omega)}+\Im(b)\|u\|_{L^2(\Omega)}^2=\langle F,\vi u\rangle_{Y^\star,Y}.
\end{gather}
Now, let us prove Theorems~\ref{thmunull} and \ref{thmbound}.
\medskip

\begin{vproof}{of Theorem~\ref{thmunull}.}
Let $F\in L^\infty(\Omega)$ with $\|F\|_{L^\infty(\Omega)}\le|a|.$ Then, the pair $(u,U)$ given by \eqref{thmunull1} is obviously a solution (Remark~\ref{rmkthmunull}). Now, assume that $(u,U)$ is a global weak solution to \eqref{nlsg} with boundary condition~\eqref{dir} or \eqref{neu}, and with $b\neq0.$ Using the dense embedding $Y\inj L^1(\Omega)$ (where $Y$ is defined as above), an appeal to \cite[Lemma~4.5]{MR3315701} (applied with $C_1=\|\nabla u\|_{L^2(\Omega)}^2+\int_\Omega\phi|u|^2\d x),$ and Hölder's inequality yield to
\begin{gather*}
\|u\|_{H^1(\Omega)}^2+\|u\|_{L^1(\Omega)}+\int_\Omega\phi|u|^2\d x\le M\|F\|_{L^\infty(\Omega)}\|u\|_{L^1(\Omega)},
\end{gather*}
for some $M=M(|a|,|b|).$ It follows that if $\|F\|_{L^\infty(\Omega)}\le\frac1M$ then
\begin{gather*}
\|u\|_{H^1(\Omega)}^2+\|u\|_{L^1(\Omega)}\le\|u\|_{L^1(\Omega)}.
\end{gather*}
Therefore, then $u=0,$ a.e.\,in $\Omega.$ From \eqref{nlsg}, we get that $a\,U=F,$ a.e.\,in $\Omega.$ Now, assume that $b=0.$ If $\Re(a)>0$ then we only deal with \eqref{predem1}. If $\Re(a)<0$ then $\Im(a)\neq0$ and so we multiply \eqref{predem2} by $\frac{|\Re(a)|+1}{\Im(a)}$ and add the result to \eqref{predem1}. In both cases, we apply Hölder's inequality and arrive at
\begin{gather*}
\|\nabla u\|_{L^2(\Omega)}^2+\|u\|_{L^1(\Omega)}+\int_\Omega\phi|u|^2\d x\le M\|F\|_{L^\infty(\Omega)}\|u\|_{L^1(\Omega)},
\end{gather*}
for some $M=M(|a|).$ The result follows by taking $\|F\|_{L^\infty(\Omega)}\le\frac1{2M}.$
\medskip
\end{vproof}

\begin{vproof}{of Theorem~\ref{thmbound}.}
Let the assumptions of the theorem be fulfilled. Let $Y$ be at the beginning of this section, and let $X=H^1_0(\Omega),$ if $u$ satisfies \eqref{dir}, and $X=Y=H^1(\Omega),$ if $u$ satisfies \eqref{neu}. By the dense embedding $Y\inj X,$ we have by \eqref{predem1}--\eqref{predem2} and \cite[Lemma~4.5]{MR3315701} that,
\begin{gather}
\label{predem3}
\|u\|_X^2+\|u\|_{L^1(\Omega)}+\int_\Omega\phi|u|^2\d x\le M_1\|F\|_{X^\star}\|u\|_X,
\end{gather}
for some $M_1=M_1(|a|,|b|).$ Applying Young's inequality, the result follows.
\medskip
\end{vproof}

\noindent
From now and until the end of this section, we assume that $(a,b)\in\C^2,$ and $\phi$ satisfies~\eqref{phi}--\eqref{pphi}. Let $\delta\in[0,1],$ and let for $n\in\N$ and $u\in L^2(\Omega),$
\begin{gather}
\label{gn}
g_n(u)=
\begin{cases}
\dfrac{u}{|u|+(n-|u|)\frac1{n^2}},		&	\text{if } |u|\le n, \medskip \\
\dfrac{u}{|u|},					&	\text{if } |u| > n,
\end{cases}
\end{gather}
\begin{gather}
\label{hn}
h_n(u)=
\begin{cases}
u,				&	\text{if } |u|\le n, \medskip \\
n\dfrac{u}{|u|},		&	\text{if } |u| > n,
\end{cases}
\end{gather}
\begin{gather}
\label{fn}
f_n=ag_n(u)+(b-\delta+\phi)h_n(u).
\end{gather}
Let $u,v\in H^1_0(\Omega).$ Since $\phi h_n(u),\phi v\in L^2(\Omega),$ we have by Remark~\ref{rmkphiu} that,
\begin{gather*}
\langle\phi h_n(u),v\rangle_{H^{-1},H^1_0}=\langle h_n(u),\phi v\rangle_{L^2,L^2}\le C\|h_n(u)\|_{L^2}\|\phi\|_{L^\infty+L^{p_\phi}}\|v\|_{H^1_0}.
\end{gather*}
We then infer that for any $u\in H^1_0(\Omega),$
\begin{gather}
\label{estphin}
\|\phi h_n(u)\|_{H^{-1}(\Omega)}\le C\|h_n(u)\|_{L^2(\Omega)}\|\phi\|_{L^\infty(\Omega)+L^{p_\phi}(\Omega)},
\end{gather}
where $C=C(N,\kappa).$ If $\Omega$ has a Lipschitz continuous boundary then the same is true with $H^1(\Omega)$ instead of $H^1_0(\Omega).$
\begin{lem}
\label{lemUu}
Let $(U_n)_{n\in\N}\subset L^\infty(\Omega),$ $U\in L^\infty(\Omega)$ and $u\in L^1_\loc(\Omega).$ We define $\omega=\big\{x\in\Omega;u(x)\neq0\big\}.$ If $U_n\underset{n\to\infty}{-\!\!\!-\!\!\!-\!\!\!-\!\!\!\weak}U \text{ in } L^\infty(\Omega)_{\w\star}$ and ${U_n}\xrightarrow[n\to\infty]{\text{a.e.\;in }\omega}\frac{u}{|u|}$ then $U=\frac{u}{|u|},$ almost everywhere in $\omega.$
\end{lem}

\begin{proof*}
Set for each $k\in\N,$ $\omega_k=\omega\cap B(0,k).$ It follows from the dominated convergence Theorem that ${U_n}_{|\omega_k}\xrightarrow[n\to\infty]{L^1(\omega_k)}\frac{u}{|u|}_{|\omega_k}.$ Let $k\in\N$ and $h\in L^\infty(\Omega)\cap L^1(\Omega)$ be defined by,
\begin{gather*}
h=
\begin{cases}
\frac{u}{|u|}-U,	&	\text{in } \omega_k,				\\
0,			&	\text{in } \Omega\setminus\omega_k.
\end{cases}
\end{gather*}
We have by Hölder's inequality that for any $n\in\N,$
\begin{align*}
	&	\; \vint_{\omega_k}\left|\frac{u}{|u|}-U\right|^2\d x		\\
   =	&	\; \Re\vint_{\omega_k}\left(\frac{u}{|u|}-U_n\right)\ovl h\,\d x+\Re\vint_\Omega\left(U_n-U\right)\ovl h\,\d x	\\
  \le	&	\; \|h\|_{L^\infty(\Omega)}\left\|U_n-\frac{u}{|u|}\right\|_{L^1(\omega_k)}+\langle U_n-U,h\rangle_{L^\infty(\Omega),L^1(\Omega)}.
\end{align*}
Passing to the limit as $n\tends\infty,$ we get that $U=\frac{u}{|u|},$ a.e.\,in $\omega_k,$ for any $k\in\N,$ from which the result follows.
\medskip
\end{proof*}

\begin{lem}
\label{lembousoldir}
Let us consider the following equation.
\begin{gather}
\label{lembousol}
-\Delta u_n+\delta u_n+f_n(u_n)=F.
\end{gather}
Let $F\in H^{-1}(\Omega)+L^\infty(\Omega).$ Assume that for each $n\in\N,$ there exists a global weak solution $u_n$ to~\eqref{lembousol} with boundary condition~\eqref{dir}. If $(u_n)_{n\in\N}$ is bounded in $H^1_0(\Omega)$ and if
\begin{gather}
\label{lembousol1}
\left(\dfrac{|u_n|^2}{|u_n|+(n-|u_n|)\frac1{n^2}}\1_{\{|u_n|\le n\}}\right)_{n\in\N} \; \text{ is bounded in } L^1(\Omega),
\end{gather}
then, up to a subsequence,
\begin{gather*}
u_n\underset{n\to\infty}{\overset{H^1_0(\Omega)_\w}{-\!\!\!-\!\!\!-\!\!\!-\!\!\!-\!\!\!-\!\!\!\weak}}u \; \text{ and } \;
g_n(u_n)\underset{n\to\infty}{\overset{L^\infty(\Omega)_{\w\star}}{-\!\!\!-\!\!\!-\!\!\!-\!\!\!-\!\!\!-\!\!\!-\!\!\!\weak}}U,
\end{gather*}
where $(u,U)$ is a global weak solution to $u$ to \eqref{nlsg} and \eqref{dir}. Finally, Symmetry Property~$\ref{sym}$ holds.
\end{lem}

\begin{proof*}
Let $(u_n)_{n\in\N}\subset H^1_0(\Omega)\cap L^1(\Omega)$ be a bounded sequence of solutions to \eqref{lembousol} and \eqref{dir}. We note that for any $n\in\N,$ $\|g_n(u_n)\|_{L^\infty(\Omega)}\le1.$ It follows that there exist a subsequence, that we still denote by $(u_n)_{n\in\N},$ $u\in H^1_0(\Omega),$ and $U\in L^\infty(\Omega)$ satisfying \eqref{defsol11} such that $u_n\underset{n\to\infty}{\overset{H^1_0(\Omega)_\w}{-\!\!\!-\!\!\!-\!\!\!-\!\!\!-\!\!\!-\!\!\!\weak}}u,$ $u_n\xrightarrow[n\to\infty]{L^2_\loc(\Omega)}u,$ $g_n(u_n)\underset{n\to\infty}{\overset{L^\infty(\Omega)_{\w\star}}{-\!\!\!-\!\!\!-\!\!\!-\!\!\!-\!\!\!-\!\!\!-\!\!\!\weak}}U,$ and $u_n\xrightarrow[n\to\infty]{\text{a.e.\;in }\Omega}u.$ By the almost pointwise convergence of $(u_n)_{n\in\N},$ \eqref{lembousol1}, Fatou's Lemma and the local compactness, we deduce that Symmetry Property~$\ref{sym}$ holds, $u\in L^1(\Omega),$ $g_n(u_n)\xrightarrow[n\to\infty]{\text{a.e.\;in }\omega}\frac{u}{|u|},$ where $\omega=\big\{x\in\Omega;u(x)\neq0\big\},$ and $h_n(u_n)\xrightarrow[n\to\infty]{L^2_\loc(\Omega)}u.$ It follows from the weak$\star$ convergence of $(g_n(u_n))_{n\in\N}$ and Lemma~\ref{lemUu} that $U$ is a saturated section associated to $u.$ Since $|h_n(u)|\le|u_n|,$ we have by \eqref{rmkphiu1} and the almost pointwise convergence of $(u_n)_{n\in\N}$ that $\phi h_n(u_n)\xrightarrow[n\to\infty]{L^2(\Omega)_\w}\phi u.$ Let $v\in H^1_0(\Omega)\cap L^1(\Omega).$ We have by \eqref{lembousol} that for any $n\in\N,$
\begin{multline}
\label{lembousolneu1}
\langle\nabla u_n,\nabla v\rangle_{L^2(\Omega),L^2(\Omega)}+\langle\delta_\star u_n,v\rangle_{L^2(\Omega),L^2(\Omega)}
			+\langle ag_n(u_n),v\rangle_{L^\infty(\Omega),L^1(\Omega)}			\\
+\langle(b-\delta_\star+\phi)h_n(u_n),v\rangle_{L^2(\Omega),L^2(\Omega)}=\langle F,v\rangle_{H^{-1}(\Omega),H^1_0(\Omega)}.
\end{multline}
We use the above convergences to pass to the limit in \eqref{lembousolneu1}. We get that $u$ satisfies \eqref{defsol21} for any $v\in  H^1_0(\Omega)\cap L^1(\Omega),$ so that $u$ is a solution to \eqref{nlsg} and \eqref{dir}. This concludes the proof.
\medskip
\end{proof*}

\noindent
A trivial adaptation of the above proof gives the following result about the homogeneous Neumann boundary condition. The details are left to the reader.

\begin{lem}
\label{lembousolneu}
Let $F\in H^1(\Omega)^\star.$ Assume that for each $n\in\N,$ there exists a global weak solution $u_n$ to~\eqref{lembousol} with boundary condition~\eqref{neu}. If $(u_n)_{n\in\N}$ is bounded in $H^1(\Omega)$ then, up to a subsequence,
\begin{gather*}
u_n\underset{n\to\infty}{\overset{H^1(\Omega)_\w}{-\!\!\!-\!\!\!-\!\!\!-\!\!\!-\!\!\!-\!\!\!\weak}}u \; \text{ and } \;
g_n(u_n)\underset{n\to\infty}{\overset{L^\infty(\Omega)_{\w\star}}{-\!\!\!-\!\!\!-\!\!\!-\!\!\!-\!\!\!-\!\!\!-\!\!\!\weak}}U,
\end{gather*}
where $(u,U)$ is a global weak solution to $u$ to \eqref{nlsg} and \eqref{neu}. Finally, Symmetry Property~$\ref{sym}$ holds.
\end{lem}

\begin{lem}
\label{lemww}
Let $n\in\N.$ Then, $f_n:H^1_0(\Omega)_\w\tends L^2(\Omega)_\w$ is weakly-weakly continuous, namely, if $u_\ell\underset{\ell\to\infty}{\overset{H^1_0(\Omega)_\w}{-\!\!\!-\!\!\!-\!\!\!-\!\!\!-\!\!\!-\!\!\!\weak}}u$ then $f_n(u_\ell)\underset{\ell\to\infty}{\overset{L^2(\Omega)_\w}{-\!\!\!-\!\!\!-\!\!\!-\!\!\!-\!\!\!-\!\!\!\weak}}f_n(u).$ If $\Omega$ has a Lipschitz continuous boundary then $f_n:H^1(\Omega)_\w\tends L^2(\Omega)_\w$ is weakly-weakly continuous.
\end{lem}

\begin{proof*}
Let $n\in\N.$ Let $u_\ell\underset{\ell\to\infty}{\overset{H^1_0(\Omega)_\w}{-\!\!\!-\!\!\!-\!\!\!-\!\!\!-\!\!\!-\!\!\!\weak}}u.$ Then, $u_\ell\xrightarrow[\ell\to\infty]{L^2_\loc(\Omega)}u,$ and $u_{\ell_k}\xrightarrow[k\to\infty]{\text{a.e.\;in }\Omega}u,$ for a subsequence $(u_{\ell_k})_{k\in\N}\subset (u_n)_{n\in\N}.$ We then have,
\begin{gather*}
f_n(u_{\ell_k})\xrightarrow[k\to\infty]{\text{a.e.\,in }\Omega}f_n(u),
\end{gather*}
and by \eqref{rmkphiu1}, $(f_n(u_\ell))_{\ell\in\N}$ is bounded in $L^2(\Omega).$ It follows that $f_n(u_{\ell_k})\underset{k\to\infty}{\overset{L^2(\Omega)_\w}{-\!\!\!-\!\!\!-\!\!\!-\!\!\!-\!\!\!-\!\!\!\weak}}f_n(u).$ The limit being necessarily $f_n(u),$ we infer that $f_n(u_\ell)\underset{\ell\to\infty}{\overset{L^2(\Omega)_\w}{-\!\!\!-\!\!\!-\!\!\!-\!\!\!-\!\!\!-\!\!\!\weak}}f_n(u),$ for the whole sequence, which is the desired result. If $\Omega$ has a Lipschitz continuous boundary then the above proof applies with $H^1(\Omega)$ in place of $H^1_0(\Omega).$
\medskip
\end{proof*}

\begin{lem}
\label{lemlaxmil}
Assume that $|\Omega|<\infty.$ Then for any $n\in\N,$ $\delta\in[0,1]$ and $F\in H^{-1}(\Omega)$ $($respectively, $\delta\in(0,1]$ and $F\in H^1(\Omega)^\star),$ there exists a global weak solution $u_n\in H^1_0(\Omega)$ to~\eqref{lembousol} with boundary condition~\eqref{dir} $($respectively, $u_n\in H^1(\Omega)$ to~\eqref{lembousol} with boundary condition\eqref{neu}$).$ Finally, Symmetry Property~$\ref{sym}$ holds.
\end{lem}

\begin{proof*}
We begin with the boundary condition \eqref{neu}. By Lax-Milgram's Theorem, we know that for any $G\in H^1(\Omega)^\star,$ there exists a unique solution $u\in H^1(\Omega)$ to $-\Delta u+\delta u=G$ which satisfies \eqref{neu}. Moreover, there exists $\alpha>0$ such that for any $G\in H^1(\Omega)^\star,$ $\left\|(-\Delta+\delta I)^{-1}G\right\|_{H^1(\Omega)}\le\alpha\|G\|_{H^1(\Omega)^\star}.$ Finally, Symmetry Property~\ref{sym} holds. Let $n\in\N.$ Set $F_n=F-f_n,$ and let us consider the following mapping $T_n$ of $H^1(\Omega)$ as follows. Set
\begin{gather*}
\begin{array}{rcccl}
T_n:H^1(\Omega)	& \xrightarrow{F_n}	& H^1(\Omega)^\star		&	\xrightarrow{(-\Delta+\delta I)^{-1}}	& H^1(\Omega),		\medskip \\
u				& \longmapsto		&		F_n(u)		&                        \longmapsto    			& (-\Delta+\delta I)^{-1}(F_n)(u).
\end{array}
\end{gather*}
Let $\rho_n=2\,\alpha\,n|\Omega|^\frac12(|a|+|b|+1+C\|\phi\|_{L^\infty+L^{p_\phi}})+\alpha\|F\|_{H^1(\Omega)^\star},$ where $C$ is given by \eqref{estphin}. It follows from \eqref{estphin} that for any $u\in H^1(\Omega),$
\begin{gather*}
\|T_n(u)\|_{H^1(\Omega)}=\left\|(-\Delta+\delta I)^{-1}(F_n)(u)\right\|_{H^1(\Omega)}\le\alpha\|F_n(u)\|_{H^1(\Omega)^\star}\le\rho_n.
\end{gather*}
Then,
\begin{gather}
\label{demlemlaxmil1}
T_n(\ovl B_{H^1(\Omega)}(0,\rho_n))\subset \ovl B_{H^1(\Omega)}(0,\rho_n),	\\
\label{demlemlaxmil2}
\ovl B_{H^1(\Omega)}(0,\rho_n) \text{  is a weakly compact subset of } H^1(\Omega)_\w.
\end{gather}
With help of Lemma~\ref{lemww}, we easily see that $T_n:H^1(\Omega)_\w\tends H^1(\Omega)_\w$ is weakly-weakly continuous. As a consequence, by \eqref{demlemlaxmil1}, \eqref{demlemlaxmil2} and a corollary of the Tychonoff fixed point Theorem (Arino, Gautier and Penot~\cite{MR0794756}, Vrabie~\cite[Theorem~1.2.11, p.6]{MR1375237}), we infer that $T_n$ admits a fixed point $u_n\in H^1(\Omega).$ The Symmetry Property~\ref{sym} is obtained by working in $\{u\in H^1(\Omega);$ $u(\vR x)=u(x),$ for a.e.\:$x\in\Omega\}$ (with the obvious modification for the odd case) in place of $H^1(\Omega).$ Working with $H^1_0(\Omega)$ instead of $H^1(\Omega),$ the boundary condition \eqref{dir} is treated in the same way, with possibly $\delta=0$ (with help of Poincaré's inequality).
\medskip
\end{proof*}

\begin{lem}[\textbf{\cite[Lemma~4.5]{MR3315701}}]
\label{lemAB}
If $(a,b)$ verifies \eqref{ab} then there exists $\delta_\star=\delta_\star(|a|,|b|)\in(0,1],$ $L=L(|a|,|b|)$ and $M=M(|a|,|b|)$ such that, if $\delta\in[0,\delta_\star]$ and $C_0,$ $C_1,$ $C_2,$ $C_3,$ $C_4$ are non-negative real numbers satisfying
\begin{gather}
\label{Re}
\big|C_1+\delta C_2+\Re(a)C_3+\big(\Re(b)-\delta\big)C_4\big|\le C_0,	\\
\label{Im}
\big|\Im(a)C_3+\Im(b)C_4\big|\le C_0,
\end{gather}
then
\begin{gather}
\label{lemestimAB}
0\le C_1+\delta C_2+LC_3+LC_4\le MC_0.
\end{gather}
\end{lem}

\begin{rmk}
\label{rmklemAB}
Actually, the conclusion in \cite[Lemma~4.5]{MR3315701} is that,
\begin{gather*}
0\le C_1+LC_3+LC_4\le MC_0.
\end{gather*}
But letting totally unchanged its proof, we obtain \eqref{lemestimAB}.
\end{rmk}

\begin{lem}
\label{lembound}
Let $(a,b)$ satisfy~\eqref{ab}. Assume that $\phi$ satisfies \eqref{phip} and let $F\in H^{-1}(\Omega).$ Then there exists $\delta_\star=\delta_\star(|a|,|b|)\in(0,1],$ $L=L(|a|,|b|)$ and $M=M(|a|,|b|)$ such that, if $\delta\in(0,\delta_\star],$  and we assume that for each $n\in\N,$ there exists a global weak solution $u_n$ of \eqref{lembousol} with boundary condition~\eqref{dir}, then we have that,
\begin{gather}
\label{lembound1}
\delta\|u_n\|_X^2+L\int_{\{|u_n|\le n\}}\dfrac{|u_n|^2}{|u_n|+(n-|u_n|)\frac1{n^2}}\d x
\le\frac{M^2}\delta\|F\|_{X^\star}^2,
\end{gather}
for any $n\in\N,$ where $X=H^1_0(\Omega).$ If $F\in H^1(\Omega)^\star$ and if for each $n\in\N,$ $u_n$ is a global weak solution to \eqref{lembousol} with boundary condition~\eqref{neu} then for any $n\in\N,$ $u_n$ satisfies \eqref{lembound1}  where $X=H^1(\Omega).$
\end{lem}

\begin{proof*}
Let $\delta_\star=\delta_\star(|a|,|b|)\in(0,1],$ $L=L(|a|,|b|)$ and $M=M(|a|,|b|)$ be given by Lemma~\ref{lemAB}. Assume that $\delta\in(0,\delta_\star].$ Let $n\in\N.$ Let $X$ be as in the lemma. Choosing $u_n$ and $\vi u_n$ as test functions in \eqref{lembound1}, we obtain
\begin{multline}
\label{demlembound1}
\|\nabla u_n\|_{L^2(\Omega)}^2+\delta\|u_n\|_{L^2(\Omega)}^2
+\Re(a)\left(\int_{\{|u_n|\le n\}}\dfrac{|u_n|^2}{|u_n|+(n-|u_n|)\frac1{n^2}}\d x+\|u_n\|_{L^1(\{|u_n|>n\})}\right)	\\
+\big(\Re(b)-\delta\big)\left(\|u_n\|_{L^2(\{|u_n|\le n\})}^2+n\|u_n\|_{L^1(\{|u_n|>n\})}\right)					\\
+\int_{\{|u_n|\le n\}}\phi|u_n|^2\d x+n\int_{\{|u_n|>n\}}\phi|u_n|\d x=\langle F,u_n\rangle_{X^\star,X},
\end{multline}
\begin{multline}
\label{demlembound2}
\Im(a)\left(\int_{\{|u_n|\le n\}}\dfrac{|u_n|^2}{|u_n|+(n-|u_n|)\frac1{n^2}}\d x+\|u_n\|_{L^1(\{|u_n|>n\})}\right)		\\
+\Im(b)\left(\|u_n\|_{L^2(\{|u_n|\le n\})}^2+n\|u_n\|_{L^1(\{|u_n|>n\})}\right)=\langle F,\vi u_n\rangle_{X^\star,X}.
\end{multline}
It follows from \eqref{demlembound1}--\eqref{demlembound2} and Lemma~\ref{lemAB} that,
\begin{gather*}
\|\nabla u_n\|_{L^2(\Omega)}^2+\delta\|u_n\|_{L^2(\Omega)}^2+L\int_{\{|u_n|\le n\}}\dfrac{|u_n|^2}{|u_n|+(n-|u_n|)\frac1{n^2}}\d x
\le M\|F\|_{X^\star}\|u_n\|_X.
\end{gather*}
Applying Young's inequality to the above, we get \eqref{lembound1}.
\medskip
\end{proof*}

\begin{lem}[\textbf{Extension}]
\label{lemexi}
Let $(\Omega_n)_{n\in\N}\subset\Omega$ be a sequence of non-decreasing open subsets of $\R^N$ such that $\cup_{n\in\N}\Omega_n=\Omega.$
\begin{enumerate}
\label{lemexi1}
\item
Let $0<p\le\infty,$ $(u_n)_{n\in\N}\subset H^1_0(\Omega_n)\cap L^p(\Omega_n),$ $(U_n)_{n\in\N}\subset L^\infty(\Omega_n)$ be a sequence of saturation sections associated to $(u_n)_{n\in\N},$ and $(\phi_n)_{n\in\N}\subset H^1_0(\Omega_n).$ Assume that there exists $C>0$ such for any $n\in\N,$
\begin{gather}
\label{lemexi11}
\|u_n\|_{H^1_0(\Omega_n)}+\|u_n\|_{L^p(\Omega_n)}+\|\nabla\phi_n\|_{L^2(\Omega_n)}\le C.
\end{gather}
Then there exist $u\in H^1_0(\Omega)\cap L^p(\Omega)$ and a saturated section $U$ associated to $u$ such that, up to subsequences $($and with no change of notation$),$
\begin{gather}
\label{lemexi12}
u_n\xrightarrow[n\to\infty]{\text{a.e.\;in }\Omega}u,	\\
\label{lemexi13}
\vlim_{n\to\infty}\langle\theta u_n,\vphi_{|\Omega_n}\rangle_{\Dr^\p(\Omega_n),\Dr(\Omega_n)}
							=\langle\theta u,\vphi\rangle_{\Dr^\p(\Omega),\Dr(\Omega)},	\\
\label{lemexi14}
\vlim_{n\to\infty}\langle U_n,\vphi_{|\Omega_n}\rangle_{\Dr^\p(\Omega_n),\Dr(\Omega_n)}=\langle U,\vphi\rangle_{\Dr^\p(\Omega),\Dr(\Omega)},	\\
\label{lemexi15}
\vlim_{n\to\infty}\langle|u_n|^2,\vphi_{|\Omega_n}\rangle_{\Dr^\p(\Omega_n),\Dr(\Omega_n)}=\langle|u|^2,\vphi\rangle_{\Dr^\p(\Omega),\Dr(\Omega)},
\end{gather}
for any $\theta$ satisfying \eqref{phi}--\eqref{pphi}, and $\vphi\in\Dr(\Omega)$. Moreover, if $\Omega=\R^N$ and $N\ge3$ then there exists $\phi\in\Dr^{1,2}(\R^N)$ such that
\begin{gather}
\label{lemexi16}
\phi_n\xrightarrow[n\to\infty]{\text{a.e.\;in }\Omega}\phi,	\\
\label{lemexi17}
\vlim_{n\to\infty}\langle\phi_n,\vphi_{|\Omega_n}\rangle_{\Dr^\p(\Omega_n),\Dr(\Omega_n)}=\langle\phi,\vphi\rangle_{\Dr^\p(\R^N),\Dr(\R^N)},	\\
\label{lemexi18}
\vlim_{n\to\infty}\langle\phi_nu_n,\vphi_{|\Omega_n}\rangle_{\Dr^\p(\Omega_n),\Dr(\Omega_n)}=\langle\phi u,\vphi\rangle_{\Dr^\p(\R^N),\Dr(\R^N)},
\end{gather}
for any $\vphi\in\Dr(\R^N).$
\label{lemexi2}
\item
Let $F\in H^{-1}(\Omega).$ Let $n\in\N.$ Let us define for any $v\in H^1_0(\Omega_n),$
\begin{gather}
\label{lemexi21}
\langle F_{|\Omega_n},v\rangle=\langle F,\wt v\rangle_{H^{-1}(\Omega),H^1_0(\Omega)},
\end{gather}
where $\wt v$ is the extension of $v$ by $0$ in $\Omega\setminus\Omega_n.$ Then for any $n\in\N,$ $F_{|\Omega_n}\in H^{-1}(\Omega_n),$ $\|F_{|\Omega_n}\|_{H^{-1}(\Omega_n)}\le\|F\|_{H^{-1}(\Omega)}$ and
\begin{gather}
\label{lemexi22}
\vlim_{n\to\infty}\langle F_{|\Omega_n},\vphi_{|\Omega_n}\rangle_{\Dr^\p(\Omega_n),\Dr(\Omega_n)}=\langle F,\vphi\rangle_{\Dr^\p(\Omega),\Dr(\Omega)},
\end{gather}
for any $\vphi\in\Dr(\Omega).$
\end{enumerate}
\end{lem}

\begin{proof*}
We begin by the first part of the lemma. Let the assumptions be fulfilled. For each $n\in\N,$ let $v_n,$ $V_n$ and $\psi_n$ be the extension by $0$ in $\Omega\setminus\Omega_n$ of $u_n,$ $U_n$ and $\phi_n,$ respectively. Let $\theta$ satisfy \eqref{phi}--\eqref{pphi}. By \eqref{lemexi11} and \eqref{rmkphiu1}, $(v_n)_{n\in\N},$ $(\theta v_n)_{n\in\N}$ and $(V_n)_{n\in\N}$ are bounded in $H^1_0(\Omega)\cap L^p(\Omega),$ in $L^2(\Omega)$ and in $L^\infty(\Omega),$ respectively. It follows that there exist $u\in H^1_0(\Omega)$ and $U\in L^\infty(\Omega)$ satisfying \eqref{defsol11}, and an extraction $(n_k)_{k\in\N}\subset(n)_{n\in\N}$ such that $v_{n_k}\underset{k\to\infty}{\overset{H^1_0(\Omega)_\w}{-\!\!\!-\!\!\!-\!\!\!-\!\!\!-\!\!\!-\!\!\!\weak}}u,$ $v_{n_k}\xrightarrow[k\to\infty]{L^2_\loc(\Omega)}u,$ $V_{n_k}\underset{k\to\infty}{\overset{L^\infty(\Omega)_{\w\star}}{-\!\!\!-\!\!\!-\!\!\!-\!\!\!-\!\!\!-\!\!\!-\!\!\!\weak}}U,$ $v_{n_k}\xrightarrow[k\to\infty]{\text{a.e.\;in }\Omega}u,$ and $\theta v_{n_k}\xrightarrow[k\to\infty]{\text{a.e.\;in }\Omega}\theta u.$ It follows that $\theta v_{n_k}\xrightarrow[k\to\infty]{L^2(\Omega)_\w}\theta u.$ Since for any $x\in\Omega,$ there exists $n_0\in\N,$ such that for any $n>n_0,$ $x\in\Omega,$ we deduce that \eqref{lemexi12} holds true, which implies with Fatou's Lemma that $u\in L^p(\Omega).$ It follows from the a.e.\,pointwise convergence of $(v_{n_k})_{k\in\N}$ that,
\begin{gather*}
V_{n_k}=\frac{v_{n_k}}{|v_{n_k}|}\xrightarrow[k\to\infty]{\text{a.e.\;in }\omega}\frac{u}{|u|}, \text{ where } \omega=\big\{x\in\Omega;u(x)\neq0\big\}.
\end{gather*}
By the weak$\star$ convergence of $(V_{n_k})_{k\in\N},$ the above limit and Lemma~\ref{lemUu}, we get that $U=\frac{u}{|u|},$ a.e.\,in $\omega,$ so that $U$ is a saturated section associated to $u.$ Let $\vphi\in\Dr(\Omega).$ Since the sets $\Omega_n$ are open, there exists $k_\star\in\N$ such that for any $k>k_\star,$ $\supp\vphi\subset\Omega_{n_k},$ so that $\vphi_{|\Omega_{n_k}}\in\Dr(\Omega_{n_k}).$ We then have for any $k>k_\star,$
\begin{gather}
\label{demlemexi1}
\langle\theta u_{n_k},\vphi_{|\Omega_{n_k}}\rangle_{\Dr^\p(\Omega_{n_k}),\Dr(\Omega_{n_k})}=\langle\theta v_{n_k},\vphi\rangle_{\Dr^\p(\Omega),\Dr(\Omega)}\xrightarrow{k\to\infty}\langle\theta u,\vphi\rangle_{\Dr^\p(\Omega),\Dr(\Omega)}.
\end{gather}
The limit \eqref{lemexi14} is obtain with the same argument. Since $(|v_n|^2)_{n\in\N}$ is bounded in $L^\frac{N}{N-2}(\Omega)$ (in $L^2(\Omega),$ if $N\le2)$ and $|v_{n_k}|^2\xrightarrow[k\to\infty]{\text{a.e.\;in }\Omega}|u|^2,$ it follows that $|v_{n_k}|^2\underset{k\to\infty}{\overset{L^\frac{N}{N-2}(\Omega)_\w}{-\!\!\!-\!\!\!-\!\!\!-\!\!\!-\!\!\!-\!\!\!-\!\!\!-\!\!\!\weak}}|u|^2$ (in $L^2(\Omega)_\w,$ if $N\le2).$ We then easily obtain \eqref{lemexi15} in the same way as in \eqref{demlemexi1}. Now, assume that $\Omega=\R^N$ and $N\ge3.$ By \eqref{lemexi11}, $(\psi_{n_k})_{k\in\N}$ is bounded in $\Dr^{1,2}(\R^N).$ It follows that there exist $\phi\in\Dr^{1,2}(\R^N)$ such that extracting another subsequence to $(n_k)_{k\in\N},$ if necessary, we have that $\psi_{n_k}\underset{k\to\infty}{\overset{L^\frac{2N}{N-2}(\R^N)_\w}{-\!\!\!-\!\!\!-\!\!\!-\!\!\!-\!\!\!-\!\!\!-\!\!\!-\!\!\!-\!\!\!\weak}}\phi,$ $\psi_{n_k}\xrightarrow[k\to\infty]{L^2_\loc(\Omega)}\phi$ and $\psi_{n_k}\xrightarrow[k\to\infty]{\text{a.e.\;in }\Omega}\phi,$ from which \eqref{lemexi16} follows. Then, \eqref{lemexi17} is obtained in the same way as in \eqref{demlemexi1}. In addition, $\psi_{n_k}v_{n_k}\xrightarrow[k\to\infty]{L^1_\loc(\Omega)}\phi u,$ from which we obtain \eqref{lemexi18}. Now, we turn out to the second property. Let $n\in\N$ and $F_{|\Omega_n}$ be defined by \eqref{lemexi21}. It is clear that $F_{|\Omega_n}\in H^{-1}(\Omega_n)$ since that if $v$ is a unitary vector of $H^1_0(\Omega_n)$ then $\wt v$ is a unitary vector of $H^1_0(\Omega),$ and
\begin{gather*}
|\langle F_{|\Omega_n},v\rangle|=|\langle F,\wt v\rangle_{H^{-1}(\Omega),H^1_0(\Omega)}|\le\|F\|_{H^{-1}(\Omega)},
\end{gather*}
so that, $\|F_{|\Omega_n}\|_{H^{-1}(\Omega_n)}\le\|F\|_{H^{-1}(\Omega)}.$ Finally, the limit \eqref{lemexi22} is obtain as in \eqref{demlemexi1}. The lemma is proved.
\medskip
\end{proof*}

\begin{vproof}{of Theorems~\ref{thmexi}.}
By Lemmas~\ref{lembousoldir}, \ref{lembousolneu}, \ref{lemlaxmil} and \ref{lembound}, it remains to show the existence of a solution satisfying the boundary condition \eqref{dir} with $|\Omega|=\infty.$ By Theorem~\ref{thmbound} and the Extension Lemma~\ref{lemexi} (applied with $\Omega_n=\Omega\cap B(0,n)),$ we obtain such a solution which satisfies \eqref{nlsg} in $\Dr^\p(\Omega).$ But the terms of the equation belong to $H^{-1}(\Omega)+L^\infty(\Omega)\inj\Dr^\p(\Omega),$ from which the result follows.
\medskip
\end{vproof}

\begin{lem}
\label{lemUu0}
Let $u_1,u_2\in L^1(\Omega),$ and let $U_1$ and $U_2$ be two saturated sections associated to $u_1$ and $u_2,$ respectively. Then $\Re\big((U_1-U_2)(\ovl{u_1-u_2})\big)\ge0,$ a.e.\,in $\Omega,$ and if
\begin{gather}
\label{lemUu01}
\Re\left(\int_\Omega(U_1-U_2)(\ovl{u_1-u_2})\d x\right)=0,
\end{gather}
then
\begin{gather}
\label{lemUu02}
(U_1-U_2)(\ovl{u_1-u_2})=0,
\end{gather}
almost everywhere in $\Omega.$
\end{lem}

\begin{rmk}
\label{rmklemUu0}
By \cite[Lemma~6.1]{MR4725781}, we already know that, $\Re\left(\int_\Omega(U_1-U_2)(\ovl{u_1-u_2})\d x\right)\ge0.$
\end{rmk}

\begin{vproof}{of Lemma~\ref{lemUu0}.}
For $j\in\{1,2\},$ let $\omega_j=\big\{x\in\Omega; \; u_j(x)\neq0\big\}.$ Let us write that,
\begin{align*}
	&	\; (U_1-U_2)(\ovl{u_1-u_2})																	\\
=	&	\; \left(U_1-\frac{u_2}{|u_2|}\right)(\ovl{-u_2})\1_{\omega_1^\co\cap\,\omega_2}
				+\left(\frac{u_1}{|u_1|}-U_2\right)\ovl{u_1}\,\1_{\omega_1\cap\,\omega_2^\co}
				+\left(\frac{u_1}{|u_1|}-\frac{u_2}{|u_2|}\right)(\ovl{u_1-u_2})\1_{\omega_1\cap\,\omega_2}				\\
=	&	\; (|u_2|-U_1\ovl{u_2})\1_{\omega_1^\co\cap\,\omega_2}+(|u_1|-U_2\ovl{u_1})\1_{\omega_1\cap\,\omega_2^\co}
				+\left(\frac{u_1}{|u_1|}-\frac{u_2}{|u_2|}\right)(\ovl{u_1-u_2})\1_{\omega_1\cap\,\omega_2}				\\
\eqdef&	\; i_1+i_2+i_3,
\end{align*}
a.e.\;in $\Omega.$ Since $|U_1\ovl{u_2}|\le|u_2|$ and $|U_2\ovl{u_1}|\le|u_1|,$ we get that $\Re(i_1)\ge0$ and $\Re(i_2)\ge0.$ We also have $\Re(i_3)\ge0$ (\cite[Corollary~5.5]{MR4340780}). After integration over $\Omega,$ we infer with \eqref{lemUu01} that, actually, $\Re(i_1)=\Re(i_2)=\Re(i_3)=0.$ Now, we see a complex number $z\in\C$ as a vector
$
\vect z=
\left(
\begin{array}{c}
\Re(z) \\
\Im(z)
\end{array}
\right)
$
of $\R^2,$ whose Euclidean norm is $|\vect z|_2=|z|.$ It follows that $\Re(i_1)=0$ may be written as,
\begin{gather*}
\vect{U_1}.\vect{u_2}=|U_1||u_2|\cos(\vect{U_1},\vect{u_2})=\Re(U_1\ovl{u_2})=|u_2|,
\end{gather*}
a.e.\:in $\omega_1^\co\cap\,\omega_2,$ where $.$ denotes the scalar product between two vectors of $\R^2.$ Since $|U_1|\le1,$ we deduce from the above that $|U_1|=\cos(\vect{U_1},\vect{u_2})=1,$ so that $U_1=\frac{u_2}{|u_2|}\1_{\omega_1^\co\cap\,\omega_2},$ and so $i_1=0.$ Arguing in the same way with $\Re(i_2)=0,$ we obtain that $i_2=0.$ Since $\Re(i_3)=0,$ let us write,
\begin{gather*}
\Re\left(\left(\frac{u_1}{|u_1|}-\frac{u_2}{|u_2|}\right)(\ovl{u_1-u_2})\right)=\frac{|u_1|+|u_2|}{|u_1||u_2|}\big(|u_1||u_2|-\Re(u_1\ovl{u_2})\big)=0.
\end{gather*}
We deduce that, $\Re(u_1\ovl{u_2})=|u_1||u_2|,$ a.e.\:in $\omega_1\cap\,\omega_2.$ This may be reformulate as,
\begin{gather*}
\vect{u_1}.\vect{u_2}=|u_1||u_2|, \text{ a.e.\:in } \omega_1\cap\,\omega_2.
\end{gather*}
We then infer that for some $k>0,$ $u_2=ku_1,$ a.e.\:in $\omega_1\cap\,\omega_2.$ Then return to $i_3,$ we see that $i_3=0.$ Hence \eqref{lemUu02}.
\medskip
\end{vproof}

\begin{vproof}{of Theorem~\ref{thmuni}.}
Let $u_1,u_2$ be as in Theorem~\ref{thmuni}. It follows that $u_1-u_2$ satisfies,
\begin{gather}
\label{demthmuni1}
-\Delta(u_1-u_2)+a(U_1-U_2)+b(u_1-u_2)+\phi(u_1-u_2)=0.
\end{gather}
We first note that if $u_1=u_2$ then by \eqref{demthmuni1} we see that $U_1=U_2.$ Choosing $a(u_1-u_2)$ as test function, we get that,
\begin{gather}
\label{demthmuni2}
\Re(a)\|\nabla u_1-\nabla u_2\|_{L^2(\Omega)}^2+|a|^2\Re(I)+\vint_\Omega(\Re(a\ovl b)+\Re(a\ovl\phi))|u_1-u_2|^2\d x=0,
\end{gather}
where $I=\int_\Omega(U_1-U_2)(\ovl{u_1-u_2})\d x.$ Since $\Re(I)\ge0$ (Remark~\ref{rmklemUu0}), we deduce from \eqref{demthmuni2} that $\Re(I)=0,$ and then $I=0$ (Lemma~\ref{lemUu0}). If $\Re(a\ovl b)+\Re(a\ovl\phi)>0,$ a.e.\,in $\Omega,$ then $u_1=u_2$ by \eqref{demthmuni2}. Otherwise, taking $u_1-u_2$ and $\vi(u_1-u_2)$ as test functions in \eqref{demthmuni1}, we obtain that,
\begin{gather*}
\|\nabla u_1-\nabla u_2\|_{L^2(\Omega)}^2+\vint_\Omega(\Re(b)+\Re(\phi))|u_1-u_2|^2\d x=0,	\\
\vint_\Omega(\Im(b)+\Im(\phi))|u_1-u_2|^2\d x=0,
\end{gather*}
since $I=0.$ Therefore $u_1=u_2.$
\medskip
\end{vproof}

\section{Inequalities for spatial localization: the proofs}
\label{secproofine}

\begin{lem}
\label{lemodi}
Let $\alpha\in(0,1],$ $\beta,\rho_0,K>0,$ and $E\in W^{1,1}(0,\rho_0;\R)$ be a non-negative solution to
\begin{gather}
\label{lemodi1}
\rho^{\beta-1}E(\rho)^{1-\alpha}\le KE^\p(\rho),
\end{gather}
for almost every $\rho\in(0,\rho_0),$ with $E(0)=0.$ Then we have
\begin{gather}
\label{lemodi2}
\forall\rho\in[0,r], \; E(\rho)=0,
\end{gather}
where $r^\beta=\left(\rho_0^\beta-K\frac\beta\alpha E(\rho_0)^\alpha\right)_+.$
\end{lem}

\begin{proof*}
We may assume that $r>0,$ otherwise there is nothing to prove. We note by \eqref{lemodi1} that $E$ is non-decreasing. Therefore, it is sufficient to prove that $E(r)=0.$ We proceed by contradiction and assume that $E(r)>0.$ Then $E>0$ over $[r,\rho_0]$ and by \eqref{lemodi1}, we have that
\begin{gather*}
\int_r^{\rho_0}\rho^{\beta-1}\d\rho\le K\int_r^{\rho_0}E^\p(\rho)E(\rho)^{\alpha-1}\d\rho,
\end{gather*}
from which we get: $\rho_0^\beta-r^\beta\le K\frac\beta\alpha(E(\rho_0)^\alpha-E(r)^\alpha).$ By definition of $r,$ this implies that $E(r)\le0,$ a contradiction.
\medskip
\end{proof*}

\begin{rmk}
\label{rmkodi}
When $\beta=0,$ the same proof gives the same result with $r=\rho_0e^{-\frac{K}{\alpha}E(\rho_0)^\alpha}.$
\end{rmk}

\begin{lem}
\label{lemodip}
Let $\alpha\in(0,1),$ $K>0$ and $\rho_1>\rho_0>0.$ Let $\eps>0$ and $E\in W^{1,1}(\rho_0,\rho_1;\R)$ be a non-negative solution to
\begin{gather}
\label{lemodip1}
E(\rho)^{1-\alpha}\le KE^\p(\rho)+\eps(\rho-\rho_0)^\frac{1-\alpha}\alpha,
\end{gather}
for almost every $\rho\in(\rho_0,\rho_1).$ Then, there exist $E_\star=E_\star(\alpha,K,\rho_0,\rho_1)$ and $\eps_\star=\eps_\star(\alpha,K)$ such that if $E(\rho_1)\le E_\star$ and $\eps\le\eps_\star$ then $E(\rho_0)=0.$
\end{lem}

\begin{proof*}
Let $E_\star=\left(\frac\alpha{2K}(\rho_1-\rho_0)\right)^\frac1\alpha$ and $\eps_\star=\frac12\left(\frac\alpha{2K}\right)^\frac{1-\alpha}\alpha.$ Let for any $\rho\in[\rho_0,\rho_1],$
\begin{gather*}
G(\rho)=\left(\frac\alpha{2K}(\rho-\rho_0)\right)^\frac1\alpha.
\end{gather*}
Then, $G\in C^1([\rho_0,\rho_1];\R),$ $G(\rho_1)=E_\star$ and $G$ satisfies,
\begin{gather*}
G(\rho)^{1-\alpha}-KG^\p(\rho)=\frac12 G(\rho)^{1-\alpha}=\eps_\star(\rho-\rho_0)^\frac{1-\alpha}\alpha,
\end{gather*}
for any $\rho\in[\rho_0,\rho_1].$ It follows from the assumption $\eps\le\eps_\star$ and \eqref{lemodip1} that,
\begin{gather}
\label{demlemodip1}
E(\rho)^{1-\alpha}-KE^\p(\rho)\le G(\rho)^{1-\alpha}-KG^\p(\rho),
\end{gather}
for almost every $\rho\in(\rho_0,\rho_1).$ Now we claim that for any $\rho\in[\rho_0,\rho_1],$ $E(\rho)\le G(\rho).$ Otherwise, by the assumption $E(\rho_1)\le G(\rho_1)$ and continuity, there would exist $r\in(\rho_0,\rho_1]$ and $\delta\in(0,r-\rho_0)$ such that $E(r)=G(r)$ and for any $\rho\in(r-\delta,r),$ $E(\rho)>G(\rho).$ This would give with \eqref{demlemodip1} that for a.e.\:$\rho\in(r-\delta,r),$ $E^\p(\rho)>G^\p(\rho).$ Integrating this expression over $(\rho,r),$ we would obtain that for any $\rho\in(r-\delta,r),$ $E(\rho)<G(\rho).$ A contradiction. Hence the claim. In particular, $E(\rho_0)\le G(\rho_0)=0.$ Hence the result.
\medskip
\end{proof*}

\begin{vproof}{of Theorems~\ref{thmsta} and \ref{thmstaF}.}
Let the assumptions of the theorems be fulfilled. Let us write $\rho_\star=\rho_0$ and $\delta=0,$ for the proof of Theorem~\ref{thmsta}, and $\rho_\star=\rho_1$ and $\delta=1,$ for the proof of Theorem~\ref{thmstaF}. Let $\rho\in(0,\rho_\star).$ We set $E(\rho)=\|\nabla u \|_{L^2(B(x_0,\rho))}^2$ and $b(\rho)=\|u\|_{L^1(B(x_0,\rho))}.$ We now proceed with the proof in 5 steps.
\\
{\bf Step~1.} $E\in W^{1,1}(0,\rho_\star),$ for a.e.\;$\rho\in(0,\rho_\star),$ $E^\p(\rho)=\|\nabla u\|_{L^2(\vsS(x_0,\rho))}^2$ and
\begin{gather}
\label{proofthmsta1}
E(\rho)+b(\rho)\le
\frac12\left(K_1(\tau)\rho^{-(N+1)}E^\p(\rho)\right)^\frac12\left(E(\rho)+b(\rho)\right)^\frac{\gamma(\tau)+1}{2}+\delta M^2\|F\|_{L^2(B(x_0,\rho))}^2,
\end{gather}
where $K_1(\tau)=C(N)M^2\max\left\{\rho_\star^{N+1},1\right\}\max\{b(\rho_\star)^{\mu(\tau)},b(\rho_\star)^{1-\gamma(\tau)}\}.$
\\
We first note that $E(\rho)=\int_0^\rho\left(\int_{\vsS(x_0,r)}|\nabla u|^2\d\sigma\right)\d r,$ for any $\rho\in(0,\rho_\star).$ The mapping
$r\longmapsto\int_{\vsS(x_0,r)}|\nabla u|^2\d\sigma$ belonging to $L^1(0,\rho_\star),$ it follows that $E$ is absolutely continuous on $(0,\rho_0).$ Then, $E\in W^{1,1}(0,\rho_\star),$ and for a.e.\;$\rho\in(0,\rho_\star),$ $E^\p(\rho)=\|\nabla u\|_{L^2(\vsS(x_0,\rho))}^2.$ Let $\rho\in(0,\rho_\star).$ By the interpolation-trace inequality (D\'iaz and Véron~\cite[Corollary~2.1]{MR792828}), we have
\begin{gather}
\label{proofthmsta2}
\|u\|_{L^2(\vsS(x_0,\rho))}\le C(N)\left(E(\rho)^\frac{1}{2}+\rho^{-\frac{N+2}2}b(\rho)\right)^\frac{N+1}{N+2}b(\rho)^\frac1{N+2},
\end{gather}
where $C=C(N).$ Applying the Cauchy-Schwarz inequality to \eqref{thmsta1} or \eqref{thmstaF1} (according to the different theorems to prove), and using \eqref{proofthmsta2}, we get that
\begin{gather}
\label{proofthmsta3}
E_T(\rho)+b(\rho)\le CME^\p(\rho)^\frac{1}{2}\left(E(\rho)^\frac{1}{2}+\rho^{-\frac{N+2}2}b(\rho)\right)^\frac{N+1}{N+2}b(\rho)^\frac1{N+2}+\delta M\vint_{B(x_0,\rho)}|F(x)u(x)|\d x,
\end{gather}
where $E_T(\rho)=E(\rho),$ for the proof of Theorem~\ref{thmsta}, and $E_T(\rho)=\|u \|_{H^1(B(x_0,\rho))}^2,$ for the proof of Theorem~\ref{thmstaF}. In the case of Theorem~\ref{thmstaF}, we apply Young's inequality to obtain
\begin{gather}
\label{proofthmsta4}
\vint_{B(x_0,\rho)}|F(x)u(x)|\d x\le\frac{M}2\|F\|_{L^2(B(x_0,\rho))}^2+\frac1{2M}\|u\|_{L^2(B(x_0,\rho))}^2.
\end{gather}
Putting together \eqref{proofthmsta3} and \eqref{proofthmsta4}, we obtain for both theorems,
\begin{gather}
\label{proofthmsta5}
E(\rho)+b(\rho)\le2CME^\p(\rho)^\frac{1}{2}\left(E(\rho)^\frac{1}{2}+\rho^{-\frac{N+2}2}b(\rho)\right)^\frac{N+1}{N+2}b(\rho)^\frac1{N+2}
+\delta M^2\|F\|_{L^2(B(x_0,\rho))}^2.
\end{gather}
Let $\tau\in\left(\frac{m+1}{2},1\right].$ A straightforward calculation yields
\begin{align*}
	&	\; \left(E(\rho)^\frac12+\rho^{-\frac{N+2}2}b(\rho)\right)b(\rho)^\frac1{N+1}							\medskip \\
  =	&	\; E(\rho)^\frac12b(\rho)^\frac1{N+1}+\rho^{-\frac{N+2}2}b(\rho)^\frac{N+2}{N+1} 					\medskip \\
  =	&	\; E(\rho)^\frac12b(\rho)^\frac\tau{N+1}b(\rho)^\frac{1-\tau}{N+1}
			+\rho^{-\frac{N+2}2}b(\rho)^{\frac12+\frac\tau{N+1}}b(\rho)^{\frac{N+2}{N+1}-\frac\tau{N+1}-\frac12}	\medskip \\
   =	&	\; E(\rho)^\frac12b(\rho)^\frac\tau{N+1}b(\rho)^{\mu(\tau)\frac{N+2}{2(N+1)}}
			+\rho^{-\frac{N+2}2}b(\rho)^{\frac12+\frac\tau{N+1}}b(\rho)^{(1-\gamma(\tau))\frac{N+2}{2(N+1)}}	\medskip \\
  \le	&	\; 2\rho^{-\frac{N+2}2}\max\left\{\rho_\star^\frac{N+2}2,1\right\}K_2^2(\tau)^\frac{N+2}{2(N+1)}
			\left(E(\rho)+b(\rho)\right)^{\frac12+\frac\tau{N+1}},
\end{align*}
where $K_2^2(\tau)=\max\{b(\rho_\star)^{\mu(\tau)},b(\rho_\star)^{1-\gamma(\tau)}\}.$ Hence~\eqref{proofthmsta1} follows from~\eqref{proofthmsta5} and the above estimate with $K_1(\tau)=64C^2M^2K_2^2(\tau)\max\left\{\rho_\star^{N+1},1\right\},$ since $\left(\frac12+\frac\tau{N+1}\right)\frac{N+1}{N+2}=\frac{\gamma(\tau)+1}{2}.$
\\
{\bf Step~2.} For any $\tau\in\left(\frac12,1\right]$ and a.e.\;$\rho\in(0,\rho_\star),$
\begin{gather*}
0\le E(\rho)^{1-\gamma(\tau)}\le K_1(\tau)\rho^{-(N+1)}E^\p(\rho)+\delta(2M)^{2(1-\gamma(\tau))}\|F\|_{L^2(B(x_0,\rho))}^{2(1-\gamma(\tau))}.
\end{gather*}
Let $\tau\in\left(\frac12,1\right]$ and $\rho\in(0,\rho_\star),$ Using the following Young inequality,
\begin{gather*}
xy\le\frac{\eps^{p^\p}}{p^\p}x^{p^\p}+\frac1{p\eps^p}y^p,
\end{gather*}
with $x=\frac12\left(K_1(\tau)\rho^{-(N+1)}E^\p(\rho)\right)^\frac12,$ $y=\left(E(\rho)+b(\rho)\right)^\frac{\gamma(\tau)+1}2,$ $p=\frac2{1+\gamma(\tau)}$ and $\eps=(\gamma(\tau)+1)^\frac{\gamma(\tau)+1}2,$ we have
\begin{align*}
	&	\; \frac12\left(K_1(\tau)\rho^{-(N+1)}E^\p(\rho)\right)^\frac12\left(E(\rho)+b(\rho)\right)^\frac{\gamma(\tau)+1}2	\\
  \le	&	\; \frac{C(\tau)}{2^\frac2{1-\gamma(\tau)}}\left(K_1(\tau)\rho^{-(N+1)}E^\p(\rho)\right)^\frac{1}{1-\gamma(\tau)}+\frac{1}{2}(E(\rho)+b(\rho)),
\end{align*}
where,
\begin{gather*}
C(\tau)=\frac{1-\gamma(\tau)}2(1+\gamma(\tau))^\frac{1+\gamma(\tau)}{1-\gamma(\tau)}\le\frac122^\frac2{1-\gamma(\tau)}.
\end{gather*}
We then obtain,
\begin{gather*}
\left(K_1(\tau)\rho^{-(N+1)}E^\p(\rho)\right)^\frac12\left(E(\rho)+b(\rho)\right)^\frac{\gamma(\tau)+1}2
\le\left(K_1(\tau)\rho^{-(N+1)}E^\p(\rho)\right)^\frac1{1-\gamma(\tau)}+(E(\rho)+b(\rho)).
\end{gather*}
Putting this estimate in \eqref{proofthmsta1}, we obtain
\begin{gather*}
E(\rho)+b(\rho)\le
\left(K_1(\tau)\rho^{-(N+1)}E^\p(\rho)\right)^\frac{1}{1-\gamma(\tau)}+2\delta M^2\|F\|_{L^2(B(x_0,\rho))}^2,
\end{gather*}
Raising both sides of this inequality to the power $(1-\gamma(\tau))\in(0,1),$ we obtain the desired result.
\\
{\bf Step~3.} Let $r\in(0,\rho_0].$ If $E(r)=0$ then $u=0,$ almost everywhere in $B(x_0,r).$
\\
It follows from the hypothesis that $E^\p=0,$ almost everywhere on $(0,r).$ We also have by definition of $\delta$ and \eqref{thmstaF2} that for both theorem, $\delta M^2\|F\|_{L^2(B(x_0,\rho))}^2=0,$ for any $\rho\in[0,r].$ It then follows from Step~1 that $b(r)=0,$ which is the desired result.
\\
{\bf Step~4.} Proof of Theorem~\ref{thmsta}.
\\
Let for $\tau\in\left(\frac12,1\right],$ $r(\tau)^{N+2}=\left(\rho_0^{N+2}-K_1(\tau)\frac{N+2}{\gamma(\tau)}E(\rho_0)^{\gamma(\tau)}\right)_+,$ where $K_1$ is given at Step~1. Let $\tau\in\left(\frac12,1\right].$ We have by Step~2 that for almost every $\rho\in(0,\rho_0),$
\begin{gather*}
\rho^{N+1}E(\rho)^{1-\gamma(\tau)}\le K_1(\tau) E^\p(\rho).
\end{gather*}
It follows from Lemma~\ref{lemodi} that for any $\rho\in[0,r(\tau)],$ $E(\rho)=0,$ and then $E(\rho_\M)=0,$ where we have set $\rho_\M=\max_{\tau\in(\frac12,1]}r(\tau).$ We conclude with Step~3.
\\
{\bf Step~5.} Proof of Theorem~\ref{thmstaF}.
\\
Let $L>0.$ Assume that $b(\rho_1)\le L.$ Set $K=K_1(1)\rho_0^{-(N+1)},$ where $K_1$ is given at Step~1. It follows that $K=K(M,L,N,\rho_1,\rho_0).$ Finally, set $\alpha=\gamma(1)=(N+2)^{-1}.$ Let then $E_\star=E_\star(M,L,N,\rho_1,\rho_0)$ and $\wt{\eps_\star}=\wt{\eps_\star}(M,L,N,\rho_1,\rho_0)$ be given by Lemma~\ref{lemodip}. Choosing $\eps_\star>0$ such that $\wt{\eps_\star}=(2M)^{2(1-\alpha)}\eps_\star^{1-\alpha},$ it follows from Step~2 and \eqref{thmstaF2} that
\begin{gather*}
E(\rho)^{1-\alpha}\le KE^\p(\rho)+\wt{\eps_\star}(\rho-\rho_0)^\frac{1-\alpha}\alpha,
\end{gather*}
for almost every $\rho\in(\rho_0,\rho_1).$ It follows from Lemma~\ref{lemodip} that if $E(\rho_1)\le E_\star$ then $E(\rho_0)=0.$ The result then comes from Step~3. This achieves the proof.
\medskip
\end{vproof}

\section{Solutions compactly supported: the proofs}
\label{seproofsolcom}

\begin{lem}
\label{lemsolcom}
Let $(a,b)$ satisfy~\eqref{ab}, let $\phi$ satisfy~\eqref{phi}--\eqref{phip}, and let $F\in L^1_\loc(\Omega).$ Let $u$ be a global weak solution to \eqref{nlsg} with the boundary condition \eqref{dir} or \eqref{neu} $($with the additional assumption for $F$ given by Definition~$\ref{defsol}).$ Let $x_0\in\Omega$ and $\rho_0>0.$ If $u$ satisfies \eqref{neu} then assume further that $\rho_0\le\dist(x_0,\Gamma).$ Then there exists $M=M(|a|,|b|)$ such that if $F_{|\Omega\cap B(x_0,\rho_0)}\in L^\infty(\Omega\cap B(x_0,\rho_0))$ with $\|F\|_{L^\infty(\Omega\cap B(x_0,\rho_0))}\le\frac1M$ then we have
\begin{gather}
\label{lemsolcom1}
\|u\|_{H^1(\Omega\cap B(x_0,\rho))}^2+\|u\|_{L^1(\Omega\cap B(x_0,\rho))}\le M\left|\dsp\int_{\Omega\cap\vsS(x_0,\rho)}u\ovl{\nabla u}.\frac{x-x_0}{|x-x_0|}\d\sigma\right|,
\end{gather}
for any $\rho\in(0,\rho_0).$
\end{lem}

\begin{proof*}
Let $U$ be the saturated section associated to the solution $u$. Let us rewrite \eqref{nlsg} as,
\begin{gather*}
-\Delta u+f(u)=G,
\end{gather*}
with $f(u)=bu,$ and $G=F-aU-\phi u.$ With help of \eqref{rmkphiu1}, we may apply \cite[Theorem~3.1]{MR3190983} to obtain,
\begin{multline*}
\|\nabla u\|_{L^2(\Omega\cap B(x_0,\rho))}^2+\Re(a)\|u\|_{L^1(\Omega\cap B(x_0,\rho))}
								+\Re(b)\|u\|_{L^2(\Omega\cap B(x_0,\rho))}^2		\\
+\int_{\Omega\cap B(x_0,\rho)}\phi|u|^2\d x
=\Re\left(\:\int_{\Omega\cap B(x_0,\rho)}F\,\ovl u\,\d x\right)+\Re\left(\:\int_{\Omega\cap\vsS(x_0,\rho)}u\ovl{\nabla u}.\frac{x-x_0}{|x-x_0|}\d\sigma\right),
\end{multline*}
\begin{multline*}
\Im(a)\|u\|_{L^1(\Omega\cap B(x_0,\rho))}+\Im(b)\|u\|_{L^2(\Omega\cap B(x_0,\rho))}^2		\\
=\Im\left(\:\int_{\Omega\cap B(x_0,\rho)}F\,\ovl u\,\d x\right)+\Im\left(\:\int_{\Omega\cap\vsS(x_0,\rho)}u\ovl{\nabla u}.\frac{x-x_0}{|x-x_0|}\d\sigma\right),
\end{multline*}
for any $\rho\in[0,\rho_0).$ Applying \cite[Lemma~4.5]{MR3315701} to the above (see also Lemma~\ref{lemAB}), and Hölder's inequality, we obtain that $u$ satisfies
\begin{multline*}
\|u\|_{H^1(\Omega\cap B(x_0,\rho))}^2+\|u\|_{L^1(\Omega\cap B(x_0,\rho))}			\\
\le C\left|\dsp\int_{\Omega\cap\vsS(x_0,\rho)}u\ovl{\nabla u}.\frac{x-x_0}{|x-x_0|}\d\sigma\right|+C\|F\|_{L^\infty(\Omega\cap B(x_0,\rho))}\|u\|_{L^1(\Omega\cap B(x_0,\rho))},
\end{multline*}
for any $\rho\in[0,\rho_0),$ where $C=C(|a|,|b|).$ The result follows by setting $M=2C.$
\medskip
\end{proof*}

\begin{vproof}{of Theorem~\ref{thmsolcomRN}.}
Let the assumption of the theorem be fulfilled. Let $M=M(|a|,|b|)$ be given by Lemma~\ref{lemsolcom} and assume that $\|F\|_{L^\infty(K^\co)}\le\frac1M.$ Let $R>0$ be such that $K\subset B(0,R).$ By Lemma~\ref{lemsolcom}, $u$ satisfies \eqref{thmsta1}, for any $B(x_0,2)\subset B(0,R)^\co.$ Choosing $\rho_0=2$ in Theorem~\ref{thmsta}, there exists $\eps_0$ such that for any $x_0\in\R^N,$ if $\|\nabla u\|_{L^2(B(x_0,2))}^2+\|u\|_{L^1(B(x_0,2))}<\eps_0$ then $\rho_\M>1,$ where $\rho_\M$ is given by \eqref{thmsta2}. Since $u\in H^1(\R^N)\cap L^1(\R^N),$ there exists $R_0>R$ such that
\begin{gather*}
\|\nabla u\|_{L^2(\{|x|>R_0\})}^2+\|u\|_{L^1(\{|x|>R_0\})}<\eps_0.
\end{gather*}
Finally, since for any $|x_0|>R_0+2,$ $B(x_0,2)\subset B(0,R_0)^\co\subset B(0,R)^\co,$ it follows from Theorem~\ref{thmsta} that for any $|x_0|>R_0+2,$ $u=0,$ a.e.\:in $B(x_0,1).$ We conclude that $u=0,$ a.e.\:in $B(0,R_0+1)^\co.$ 
\medskip
\end{vproof}

\begin{vproof}{of Theorems~\ref{thmsolcom}.}
Let $M_1=M_1(|a|,|b|)$ and $M_2=M_2(|a|,|b|)$ be given by \eqref{predem3} and Lemma~\ref{lemsolcom}, respectively. Let $M=M_1+M_2,$ let $K\subset\Omega$ be a compact subset, $F\in H^{-1}(\Omega),$ (respectively, $F\in H^1(\Omega)^\star)$ with $\|F\|_{L^\infty(\Omega\setminus K)}\le\frac1M,$ and let $u$ be a global weak solution to \eqref{nlsg} satisfying one of the two boundary conditions \eqref{dir} or \eqref{neu}. Finally, let $\eps_\star>0$ be such that $K(5\eps_\star)\subset\Omega.$ Let $\eps\in(0,\eps_\star).$ By Lemma~\ref{lemsolcom}, we have that $u$ satisfies \eqref{thmsta1}, for almost every $\rho\in(0,2\eps),$ and any $x_0\in\Omega$ such that $K\cap B(x_0,2\eps)=\emptyset$ and $B(x_0,2\eps)\subset\Omega.$ By Theorem~\ref{thmbound}, there exists $\delta=\delta(\eps,|a|,|b|,M)$ such that if $\|F\|_{H^1(\Omega)^\star}<\delta$ then $\rho_\M>\eps,$ where $\rho_\M$ is given by \eqref{thmsta2}. We have that,
\begin{gather}
\label{demthmsolcom1}
K\cap B(x_0,2\eps)=\emptyset, \text{ for any } x_0\in\Omega\setminus K(2\eps),
\end{gather}
and $B(x_0,2\eps)\subset\Omega,$ for any $x_0\in K(3\eps).$ We may use Theorem~\ref{thmsta} to conclude that $u=0,$ a.e.\:in $B(x_0,\eps),$ for any $x_0\in K(3\eps)\setminus K(2\eps).$ It follows that
\begin{gather}
\label{demthmsolcom2}
u=0, \text{ almost everywhere in } K(4\eps)\setminus K(\eps).
\end{gather}
Since $K(2\eps)\cap\ovl{K(3\eps)^\co}=\emptyset,$ it follows from \eqref{demthmsolcom2} that we may define $\wt u\in H^1_0(\Omega)\cap L^1(\Omega),$ if $u$ satisfies \eqref{dir}, and $\wt u\in H^1(\Omega),$ if $u$ satisfies \eqref{neu}, as
\begin{gather*}
\wt u=
\begin{cases}
u,	&	\text{in } \Omega_\eps\eqdef\Omega\setminus K(3\eps),	\\
0,	&	\text{in } K(3\eps).
\end{cases}
\end{gather*}
Choosing $\wt u$ and $\vi\wt u$ as test functions, we obtain  that
\begin{gather}
\label{demthmsolcom3}
\|\nabla u\|_{L^2(\Omega_\eps)}^2+\Re(a)\|u\|_{L^1(\Omega_\eps)}+\Re(b)\|u\|_{L^2(\Omega_\eps)}^2+\int_{\Omega_\eps}\phi|u|^2\d x
\le\int_{\Omega_\eps}|F\,\ovl u|\d x,	\\
\label{demthmsolcom4}
\Im(a)\|u\|_{L^1(\Omega_\eps)}+\Im(b)\|u\|_{L^2(\Omega_\eps)}^2\le\int_{\Omega_\eps}|F\,\ovl u|\d x.
\end{gather}
It follows from \cite[Lemma~4.5]{MR3315701} that there exists $M_3=M_3(|a|,|b|)$ such that
\begin{gather}
\label{demthmsolcom5}
\|u\|_{H^1(\Omega_\eps)}^2+\|u\|_{L^1(\Omega_\eps)}+\int_{\Omega_\eps}\phi|u|^2\d x\le M_3\int_{\Omega_\eps}|F\,\ovl u|\d x.
\end{gather}
Now, let us note that the constant $M_1$ in \eqref{predem3} is obtain from \cite[Lemma~4.5]{MR3315701} applied to \eqref{predem1}--\eqref{predem2}, in which the constants involved $(\Re(a),$ $\Im(a),$ $\Re(b),$ $\Im(b)$ and $\delta=0)$ are exactly the same as in \eqref{demthmsolcom3}--\eqref{demthmsolcom4}. We then infer that $M_3=M_1.$ Applying Hölder's inequality, it follows from \eqref{demthmsolcom5} that,
\begin{gather*}
\|u\|_{H^1(\Omega_\eps)}^2+(1-M_1\|F\|_{L^\infty(\Omega\setminus K)})\|u\|_{L^1(\Omega_\eps)}\le0.
\end{gather*}
By assumption, $\|F\|_{L^\infty(\Omega\setminus K)}\le\frac{1}{M_1}.$ Therefore, $u=0$ in $\Omega\setminus K(3\eps),$ hence in $\Omega\setminus K(\eps)$ by \eqref{demthmsolcom2}. The case in which $u$ satisfies \eqref{dir} is obtained in the same way and the details are left to the reader.
\medskip
\end{vproof}

\section{Application to a Schrödinger-Poisson system}
\label{SP}

\begin{thm}[\textbf{Existence, a priori bound and compactness}]
\label{thmSPRN}
Let $(a,b)$ satisfy \eqref{ab} and let $e\ge0.$ Assume that $N\in\{3,4\}.$ Then for any $F\in H^{-1}(\R^N),$ there exists a global weak solution
$(u,U,\phi)$ to \eqref{bcRN}--\eqref{nlsSPRN}. In addition, there exists  $M=M(|a|,|b|)$ such that any global weak solution $(u,U,\phi)$ to \eqref{bcRN}--\eqref{nlsSPRN} satisfies the following properties$:$
\begin{gather}
\label{thmSPboundRN1}
\|u\|_{H^1(\R^N)}^2+\|u\|_{L^1(\R^N)}+e\int_{\R^N}\phi|u|^2\d x\le M\|F\|_{H^{-1}(\R^N)}^2,	\\
\label{thmSPboundRN2}
\|\nabla\phi\|_{L^2(\R^N)}^2=\frac{e}2\int_{\R^N}\phi|u|^2\d x.
\end{gather}
Finally, if there exists a compact subset $K$ of $\R^N$ such that $F_{|K^\co}\in L^\infty(K^\co)$ and $\|F\|_{L^\infty(K^\co)}\le\frac1M$ then $\supp u$ is compact.
\end{thm}

\begin{rmk}
\label{rmkthmSPRN}
We do not know if the solution $(u,U,\phi)$ is unique. On the other hand, if $-\Delta\phi_1=-\Delta\phi_2=\frac{e}2|u|^2,$ in $L^2(\R^N),$ then we easily obtain that $\nabla\phi_1=\nabla\phi_2,$ in $L^2(\R^N).$ Since $\phi_1,\phi_2\in L^\frac{2N}{N-2}(\R^N),$ we infer that $\phi_1=\phi_2.$ It follows that uniqueness of $u$ implies uniqueness of $\phi$ and then $U.$ In addition, the solution to $-\Delta\phi=\frac{e}2|u|^2,$ in $\Dr^\p(\R^N),$ is given by,
\begin{gather*}
\phi=\frac{e}2(-\Delta)^{-1}|u|^2=\frac{e}2G\star|u|^2\in L^1_\loc(\R^N),	\\
G(x)=\frac1{N(N-2)|B(0,1)|}\frac1{|x|^{N-2}}, \; x\neq0.
\end{gather*}
In particular, $\phi\ge0$ in $\R^N.$ By interior elliptic regularity (Cazenave~\cite[Proposition~4.1.2]{caz-sle}), we easily obtain that $\phi\in H^2_\loc(\R^N;\R).$
\end{rmk}

\begin{prop}
\label{propSP}
Assume $|\Omega|<\infty.$ Let $(a,b)$ satisfy \eqref{ab} and let $e\ge0.$ If $N\le4$ then for any $F\in H^{-1}(\Omega),$ there exists a global weak solution
\begin{gather}
\label{bc}
\begin{cases}
(u,U,\phi)\in H^1_0(\Omega)\times L^\infty(\Omega)\times H^1_0(\Omega;\R),	\medskip \\
U \text{ is a saturated section associated to } u,							\medskip \\
\phi u\in L^2(\Omega) \text{ and } \phi\ge0 \text{ in  } \Omega,
\end{cases}
\end{gather}
to
\begin{gather}
\label{nlsSP}
\begin{cases}
-\Delta u+a\,U+b\,u+e\,\phi\,u=F, \text{ in } H^{-1}(\Omega),	\medskip \\
-\Delta\phi=\dfrac{e}2|u|^2, \text{ in } L^2(\Omega).
\end{cases}
\end{gather}
In addition, any global weak solution $(u,U,\phi)$ to \eqref{bc}--\eqref{nlsSP} satisfies that
\begin{gather}
\label{propSPbound1}
\|u\|_{H^1_0(\Omega)}^2+\|u\|_{L^1(\Omega)}+e\int_\Omega\phi|u|^2\d x\le M\|F\|_{H^{-1}(\Omega)}^2,	\\
\label{propSPbound2}
\|\nabla\phi\|_{L^2(\Omega)}^2=\frac{e}2\int_\Omega\phi|u|^2\d x,
\end{gather}
for some $M=M(|a|,|b|).$
\end{prop}

\begin{proof*}
Since $N\le4$ we note that if $v\in H^1_0(\Omega)$ then by Sobolev' embedding, $v$ satisfies \eqref{phi}--\eqref{pphi} and $|v|^2\in L^2(\Omega).$ Let $F\in H^{-1}(\Omega).$ Let $\phi_0=0$ and let $(u_1,U_1)\in H^1_0(\Omega)\times L^\infty(\Omega)$ be a global weak solution to,
\begin{gather*}
-\Delta u_1+ a\,U_1+b\,u_1+e\,\phi_0u_1=F, \text{ in } H^{-1}(\Omega),
\end{gather*}
given by Theorem~\ref{thmexi}. Now, let then $\phi_1\in H^1_0(\Omega;\R)$ be a solution to $-\Delta\phi_1=\frac{e}2|u_1|^2.$ By the weak maximum principle, $\phi_1$ is non-negative. Then $\phi_1$ satisfies \eqref{phi}--\eqref{pphi}, and by Theorem~\ref{thmexi} there exists a global weak solution $(u_2,U_2)\in H^1_0(\Omega)\times L^\infty(\Omega)$ to,
\begin{gather*}
-\Delta u_2+a\,U_2+b\,u_2+e\,\phi_1u_2=F, \text{ in } H^{-1}(\Omega).
\end{gather*}
And now, we may find $\phi_2\in H^1_0(\Omega;\R)$ a non-negative solution to $-\Delta\phi_2=\frac{e}2|u_2|^2.$ By induction, we construct a sequence $(u_n,U_n,\phi_{n-1})_{n\in\N}$ of global weak solutions to \eqref{bc} and
\begin{gather}
\label{demnlsSP}
\begin{cases}
-\Delta u_n+a\,U_n+b\,u_n+e\,\phi_{n-1}u_n=F, \text{ in } H^{-1}(\Omega),	\medskip \\
-\Delta\phi_n=\dfrac{e}2|u_n|^2, \text{ in } L^2(\Omega),
\end{cases}
\end{gather}
for any $n\in\N.$ By \eqref{thmbound11}, $(u_n)_{n\in\N}$ is bounded in $H^1_0(\Omega).$ We infer from the second equation in \eqref{demnlsSP}, Cauchy-Schwarz' inequality, the Sobolev embedding $H^1_0(\Omega)\inj L^4(\Omega),$ Poincaré's inequality, and Young's inequality that,
\begin{gather*}
\|\nabla\phi_n\|_{L^2(\Omega)}^2\le\frac{e}2\|u_n\|_{L^4(\Omega)}^2\|\phi_n\|_{L^2(\Omega)}\le C+\frac12\|\nabla\phi_n\|_{L^2(\Omega)}^2.
\end{gather*}
As a consequence, $(\phi_n)_{n\in\N}$ is bounded in $H^1_0(\Omega),$ and by Sobolev' embedding, $(\phi_{n-1}u_n)_{n\in\N}$ is bounded in $L^2(\Omega).$  Then, up to a subsequence, there exist $u\in H^1_0(\Omega),$ $U\in L^\infty(\Omega)$ satisfying \eqref{defsol11}, and $\phi\in H^1_0(\Omega;\R)$ with $\phi\ge0$ such that, up to a subsequence, $u_n\underset{n\to\infty}{\overset{H^1_0(\Omega)_\w}{-\!\!\!-\!\!\!-\!\!\!-\!\!\!-\!\!\!-\!\!\!\weak}}u,$ $u_n\xrightarrow[n\to\infty]{L^2_\loc(\Omega)}u,$ $U_n\underset{n\to\infty}{\overset{L^\infty(\Omega)_{\w\star}}{-\!\!\!-\!\!\!-\!\!\!-\!\!\!-\!\!\!-\!\!\!-\!\!\!\weak}}U,$ $u_n\xrightarrow[n\to\infty]{\text{a.e.\;in }\Omega}u,$ $\phi_n\xrightarrow[n\to\infty]{L^2_\loc(\Omega)}\phi,$ and $\phi_n\xrightarrow[n\to\infty]{\text{a.e.\;in }\Omega}\phi.$ Since $|u_n|^2\xrightarrow[n\to\infty]{L^1_\loc(\Omega)}|u|^2,$ we may pass to the limit in the second equation of \eqref{demnlsSP} in $\Dr^\p(\Omega)$ to get, $-\Delta\phi=\dfrac{e}2|u|^2,$ in $L^2(\Omega).$ In addition, we have $U_n=\frac{u_n}{|u_n|}\xrightarrow{n\to\infty}\frac{u}{|u|},$ a.e.\,where $u\neq0.$ We get with help of Lemma~\ref{lemUu} that $U$ is a saturated section associated to $u.$ Finally, $(\phi_{n-1}u_n)_{n\in\N}$ is bounded in $L^2(\Omega)$ and $\phi_{n-1}u_n\xrightarrow[n\to\infty]{\text{a.e.\;in }\Omega}\phi u,$ so that $\phi u\in L^2(\Omega)$ and $\phi_{n-1}u_n\underset{n\to\infty}{\overset{L^2(\Omega)_\w}{-\!\!\!-\!\!\!-\!\!\!-\!\!\!-\!\!\!\weak}}\phi u.$ Finally, we use all these converges to pass to the limit in the first equation in \eqref{demnlsSP} in $\Dr^\p(\Omega).$ It follows that $(u,U,\phi)$ is a solution to \eqref{bc}--\eqref{nlsSP}. Taking the $H^{-1}-H^1_0$ duality product of the second equation in \eqref{nlsSP} with $\phi,$ we get \eqref{propSPbound2}. To conlude, we note that \eqref{propSPbound1} comes from \eqref{thmbound11}.
\medskip
\end{proof*}

\begin{prop}
\label{propSPcom}
Let $N\le4$ and assume $|\Omega|<\infty.$ Let $(a,b)$ satisfy \eqref{ab}, let $K\subset\Omega$ be any compact subset of $\R^N,$ let $F\in H^{-1}(\Omega)$ be such that $F_{|\Omega\setminus K}\in L^\infty(\Omega\setminus K),$ and let $(u,U,\phi)$ be any global weak solution to~\eqref{bc}--\eqref{nlsSP}. Then, there exist $M=M(|a|,|b|)$ and $\eps_\star=\eps_\star(\dist(K,\Gamma))$ such that, for any $\eps\in(0,\eps_\star),$ there exists $\delta=\delta(\eps,|a|,|b|,N)$ verifying that if $\|F\|_{H^{-1}(\Omega)}\le\delta$ and $\|F\|_{L^\infty(\Omega\setminus K)}\le\frac1M$ then $\supp u\subset K(\eps)\subset\Omega,$ where $K(\eps)$ is given by \eqref{K}.
\end{prop}

\begin{proof*}
Apply Theorem~\ref{thmsolcom}.
\medskip
\end{proof*}

\begin{vproof}{of Theorem~\ref{thmSPRN}.}
We first note that if $N\in\{3,4\}$ then $L^1(\R^N)\cap L^\frac{N}{N-1}(\R^N)\inj L^\frac{2N}{N+2}(\R^N)$ with dense embedding so that by duality,
\begin{gather}
\label{demthmSPRN}
\Dr^{1,2}(\R^N)\inj L^\frac{2N}{N-2}(\R^N)\inj L^\infty(\R^N)+L^{p_\phi}(\R^N).
\end{gather}
With help of \eqref{rmkphiu1}, it follows that if $(u,\phi)\in H^1(\R^N)\times\Dr^{1,2}(\R^N)$ then $\phi$ satisfies \eqref{phi}--\eqref{pphi} and $\phi u\in L^2(\R^N).$ Moreover, we have by \eqref{thmbound11} that if $(u,U,\phi)$ is a global weak solution to \eqref{bcRN}--\eqref{nlsSPRN} then $(u,\phi)$ satisfies \eqref{thmSPboundRN1}. Taking the $H^{-1}-H^1$ duality product of the second equation in \eqref{nlsSPRN} with $\phi,$ we obtain \eqref{thmSPboundRN2}. Finally, the compactness Property comes from Theorem~\ref{thmsolcomRN}. It remains to show the existence of a solution. Let $F\in H^{-1}(\R^N).$ Let $(F_{|\Omega_n})_{n\in\N}\subset H^{-1}(\Omega_n)$ be defined as in the second part of the Extension Lemma~\ref{lemexi}, where $\Omega_n=B(0,n),$ and let us apply Proposition~\ref{propSP}. For each $n\in\N,$ let $(u_n,U_n,\phi_n)$ be a global weak solution to \eqref{bc}--\eqref{nlsSP}, where the domain is $\Omega_n$ and where the right member of the first equation in \eqref{nlsSP} is $F_{|\Omega_n}.$ By \eqref{propSPbound1}--\eqref{propSPbound2}, and the second part of the Extension Lemma~\ref{lemexi}, we have that 
\begin{gather*}
\|u_n\|_{H^1_0(\Omega_n)}^2+\|u_n\|_{L^1(\Omega_n)}+\|\nabla\phi_n\|_{L^2(\Omega_n)}^2\le M\|F_{|\Omega_n}\|_{H^{-1}(\Omega_n)}^2\le M\|F\|_{H^{-1}(\R^N)}^2,
\end{gather*}
for any $n\in\N.$ Therefore, we may apply the first part of the Extension Lemma~\ref{lemexi} to obtain the existence of a global weak solution $(u,U,\phi)$ to \eqref{bcRN}--\eqref{nlsSPRN} in $\Dr^\p(\R^N).$ It is clear that the second equation in \eqref{nlsSPRN} makes sense in $L^2(\R^N),$ while all the terms of the first equation belong to $H^{-1}(\R^N)+L^\infty(\R^N)\inj\Dr^\p(\R^N).$ This ends the proof of the theorem.
\medskip
\end{vproof}

\baselineskip .4cm


\begin{thebibliography}{10}
\addcontentsline{toc}{section}{References}

\bibitem{ak}
G.~P. Agrawal and Y.~S. Kivshar.
\newblock {\em Optical {S}olitons: {F}rom {F}ibers to {P}hotonic {C}rystals}.
\newblock Academic Press, California, San Diego, 2003.

\bibitem{MR2447960}
A.~Ambrosetti.
\newblock Remarks on some systems of nonlinear {S}chr\"odinger equations.
\newblock {\em J. Fixed Point Theory Appl.}, 4(1):35--46, 2008.

\bibitem{MR2002i:35001}
S.~N. Antontsev, J.~I. D{\'{\i}}az, and S.~Shmarev.
\newblock {\em Energy methods for free boundary problems}.
\newblock Progress in Nonlinear Differential Equations and their Applications,
  48. Birkh\"auser Boston Inc., Boston, MA, 2002.
\newblock Applications to nonlinear PDEs and fluid mechanics.

\bibitem{MR0794756}
O.~Arino, S.~Gautier, and J.-P. Penot.
\newblock A fixed point theorem for sequentially continuous mappings with
  application to ordinary differential equations.
\newblock {\em Funkcial. Ekvac.}, 27(3):273--279, 1984.

\bibitem{MR4521439}
P.~B\'{e}gout.
\newblock The dual space of a complex {B}anach space restricted to the field of
  real numbers.
\newblock {\em Adv. Math. Sci. Appl.}, 31(2):241--252, 2022.

\bibitem{MR2876246}
P.~B{\'e}gout and J.~I. D{\'{\i}}az.
\newblock Localizing estimates of the support of solutions of some nonlinear
  {S}chr\"odinger equations --- {T}he stationary case.
\newblock {\em Ann. Inst. H. Poincar\'e Anal. Non Lin\'eaire}, 29(1):35--58,
  2012.

\bibitem{MR3190983}
P.~B{\'e}gout and J.~I. D{\'{\i}}az.
\newblock A sharper energy method for the localization of the support to some
  stationary {S}chr\"odinger equations with a singular nonlinearity.
\newblock {\em Discrete Contin. Dyn. Syst.}, 34(9):3371--3382, 2014.

\bibitem{MR3315701}
P.~B{\'e}gout and J.~I. D{\'{\i}}az.
\newblock Existence of weak solutions to some stationary {S}chr\"odinger
  equations with singular nonlinearity.
\newblock {\em Rev. R. Acad. Cienc. Exactas F\'\i s. Nat. Ser. A Math. RACSAM},
  109(1):43--63, 2015.

\bibitem{MR4340780}
P.~B\'{e}gout and J.~I. D\'{\i}az.
\newblock Finite time extinction for a class of damped {S}chr\"{o}dinger
  equations with a singular saturated nonlinearity.
\newblock {\em J. Differential Equations}, 308:252--285, 2022.

\bibitem{MR4725781}
P.~B\'{e}gout and J.~I. D\'{\i}az.
\newblock Strong stabilization of damped nonlinear {S}chr\"{o}dinger equation
  with saturation on unbounded domains.
\newblock {\em J. Math. Anal. Appl.}, 538(1):Paper No. 128329, 2024.

\bibitem{MR1659454}
V.~Benci and D.~Fortunato.
\newblock An eigenvalue problem for the {S}chr\"{o}dinger-{M}axwell equations.
\newblock {\em Topol. Methods Nonlinear Anal.}, 11(2):283--293, 1998.

\bibitem{MR0482275}
J.~Bergh and J.~L{\"o}fstr{\"o}m.
\newblock {\em Interpolation spaces. {A}n introduction}.
\newblock Springer-Verlag, Berlin, 1976.
\newblock Grundlehren der Mathematischen Wissenschaften, No. 223.

\bibitem{MR2272971}
A.~Biswas and S.~Konar.
\newblock {\em Introduction to non-{K}err law optical solitons}.
\newblock Chapman \& Hall/CRC Applied Mathematics and Nonlinear Science Series.
  Chapman \& Hall/CRC, Boca Raton, FL, 2007.

\bibitem{MR0481460}
H.~Brezis.
\newblock Solutions of variational inequalities, with compact support.
\newblock {\em Uspehi Mat. Nauk}, 29(2(176)):103--108, 1974.
\newblock Translated from the English by Ju. A. Dubinski\u i.

\bibitem{MR2759829}
H.~Brezis.
\newblock {\em Functional analysis, {S}obolev spaces and partial differential
  equations}.
\newblock Universitext. Springer, New York, 2011.

\bibitem{caz-sle}
T.~Cazenave.
\newblock {\em An introduction to semilinear elliptic equations}.
\newblock Editora do Instituto de Matem{\'a}tica, Universidade Federal do Rio
  de Janeiro, Rio de Janeiro, 2006.

\bibitem{csbk}
Y.~Choukroun, A.~Shtern, A.~Bronstein, and R.~Kimmel.
\newblock Hamiltonian operator for spectral shape analysis.
\newblock 26(2):1320--1331, 2018.

\bibitem{MR3639295}
M.~Colin and T.~Watanabe.
\newblock Standing waves for the nonlinear {S}chr\"{o}dinger equation coupled
  with the {M}axwell equation.
\newblock {\em Nonlinearity}, 30(5):1920--1947, 2017.

\bibitem{MR853732}
J.~I. D{\'{\i}}az.
\newblock {\em Nonlinear partial differential equations and free boundaries.
  {V}ol. {I}}, volume 106 of {\em Research Notes in Mathematics}.
\newblock Pitman (Advanced Publishing Program), Boston, MA, 1985.
\newblock Elliptic equations.

\bibitem{MR2466410}
J.~I. D{\'{\i}}az.
\newblock Estimates of the location of a free boundary for the obstacle and
  {S}tefan problems obtained by means of some energy methods.
\newblock {\em Georgian Math. J.}, 15(3):475--484, 2008.

\bibitem{MR792828}
J.~I. D{\'{\i}}az and L.~V{\'e}ron.
\newblock Local vanishing properties of solutions of elliptic and parabolic
  quasilinear equations.
\newblock {\em Trans. Amer. Math. Soc.}, 290(2):787--814, 1985.

\bibitem{MR2808162}
G.~P. Galdi.
\newblock {\em An introduction to the mathematical theory of the
  {N}avier-{S}tokes equations}.
\newblock Springer Monographs in Mathematics. Springer, New York, second
  edition, 2011.
\newblock Steady-state problems.

\bibitem{MR1817225}
E.~H. Lieb and M.~Loss.
\newblock {\em Analysis}, volume~14 of {\em Graduate Studies in Mathematics}.
\newblock American Mathematical Society, Providence, RI, second edition, 2001.

\bibitem{MR3016511}
L.~A. Maia, E.~Montefusco, and B.~Pellacci.
\newblock Weakly coupled nonlinear {S}chr\"odinger systems: the saturation
  effect.
\newblock {\em Calc. Var. Partial Differential Equations}, 46(1-2):325--351,
  2013.

\bibitem{olco}
V.~Ozolins, R.~Lai, R.~Caflisch, and S.~Osher.
\newblock Compressed modes for variational problems in mathematics and physic.
\newblock 110:18368--18373, 2013.

\bibitem{MR306715}
W.~A. Strauss.
\newblock On weak solutions of semi-linear hyperbolic equations.
\newblock {\em An. Acad. Brasil. Ci.}, 42:645--651, 1970.

\bibitem{MR2296978}
F.~Tr{\`e}ves.
\newblock {\em Topological vector spaces, distributions and kernels}.
\newblock Dover Publications Inc., Mineola, NY, 2006.
\newblock Unabridged republication of the 1967 original.

\bibitem{MR1375237}
I.~I. Vrabie.
\newblock {\em Compactness methods for nonlinear evolutions}, volume~75 of {\em
  Pitman Monographs and Surveys in Pure and Applied Mathematics}.
\newblock Longman Scientific \& Technical, Harlow; copublished in the United
  States with John Wiley \& Sons, Inc., New York, second edition, 1995.
\newblock With a foreword by A. Pazy.

\end{thebibliography}

\end{document}